\documentclass[reqno]{amsart}

\usepackage[a4paper,margin=1.5in]{geometry}

\pdfoutput=1

\usepackage{mathrsfs}
\usepackage{amsfonts}
\usepackage{pdfpages}
\usepackage{amsmath}
\usepackage{amssymb}
\usepackage[foot]{amsaddr}

\usepackage{mathtools}
\usepackage[inline]{enumitem}
\usepackage{hyperref}

\newcommand{\pres}[3]{\textnormal{#1} \langle #2 \mid #3 \rangle}

\newcommand{\xra}[1]{\xrightarrow{}^\ast_{#1}}
\newcommand{\xr}[1]{\xrightarrow{}_{#1}}

\newcommand{\fC}{\mathfrak{C}}

\newcommand{\N}{\mathbb{N}}
\newcommand{\Z}{\mathbb{Z}}

\newcommand{\bS}{\mathbf{S}}
\newcommand{\dbot}{\mathbin{\text{$\bot\mkern-8mu\bot$}}}

\newtheorem{innercustomgeneric}{\customgenericname}
\providecommand{\customgenericname}{}
\newcommand{\newcustomtheorem}[2]{%
  \newenvironment{#1}[1]
  {%
   \renewcommand\customgenericname{#2}%
   \renewcommand\theinnercustomgeneric{##1}%
   \innercustomgeneric
  }
  {\endinnercustomgeneric}
}
\newtheorem*{theorem*}{Theorem} % 1st argument is your name for it
\newcustomtheorem{customthm}{Theorem}
\newcustomtheorem{customcor}{Corollary}
\newcustomtheorem{customprop}{Proposition}

\numberwithin{lemma}{section}
\numberwithin{proposition}{section}

\newtheorem*{mainlemma*}{Main Lemma}

\theoremstyle{definition}

\newtheorem*{question*}{Question}
\newtheorem{example}{Example}
\numberwithin{example}{section}

\begin{document}

\title[Tseytin's seven-relation semigroup]{G.\ S.\ Tseytin's seven-relation semigroup with undecidable word problem}

\author{\small{Carl-Fredrik Nyberg-Brodda}}
\address{{School of Mathematics, Korea Institute for Advanced Study (KIAS), Seoul 02455, Korea}}
%\curraddr{Department of Mathematics and Statistics,
%Case Western Reserve University, Cleveland, Ohio 43403}
\email{cfnb@kias.re.kr}
\thanks{The author was supported by Mid-Career Researcher Program (RS-2023-00278510) through the National Research Foundation funded by the government of Korea, and is currently supported by the KIAS Individual Grant (MG094701) at Korea Institute for Advanced Study.}

%    General info
\subjclass[2020]{}

\date{\today}

%\dedicatory{Dedicated to the memory of G. S. Tseytin (1936--2022)}

%\keywords{Word problem; semigroups; undecidability; algorithmic problems}

\begin{abstract}
We give an introduction to the ideas behind G.\ S.\ Tseytin's 1958 construction of a seven-relation semigroup with undecidable word problem. We give a history of the ideas leading up to its construction, some intuition for the proof, and provide an overview of some subsequent results and developments which stem from this remarkable semigroup. An English translation of Tseytin's article \cite{Tseytin1958} by the author of the present article is supplemented at the end of the article (the Russian original was published in \textit{Trudy Mat. Inst. Steklov} \textbf{52} (1958), pp. 172--189).
\end{abstract}

\maketitle

\noindent Can computers do algebra? More specifically, can a computer faithfully represent and compute with algebraic structures which are finitely definable? These questions lie at the heart of the \textit{word problem}, one of the oldest and most remarkable algorithmic problems in mathematics. It arises from the general problem of asking whether algebraic structures and their elements can be encoded into a computer in such a way as to make their operations computable. The first question one must ask in this line is then the following: given two representations of elements in the encoding, can we decide whether the two represented elements are equal? This is the word problem. This fascinating problem has attracted some of the greatest minds in theoretical computer science and group theory of the 20th century, including M.\ Dehn \cite{Dehn1911}, M.\ Hall Jr.\ \cite{Hall1949}, D.\ Knuth \cite{Knuth1970}, W.\ Magnus \cite{Magnus1932}, G.\ S.\ Makanin \cite{Makanin1966}, Yu.\ V.\ Matiyasevich \cite{Matiyasevich1967}, P.\ S.\ Novikov \cite{Novikov1952}, A.\ Turing \cite{Turing1950}, E.\ L.\ Post \cite{Post1947}, D.\ Scott \cite{Scott1956}, and A.\ Thue \cite{Thue1914}, all of whom worked directly with the word problem for groups or semigroups. It should come as no surprise, then, that G.\ S.\ Tseytin also worked on this problem. 

The long list of names of people who worked on the word problem above is not merely an exercise in literature review, but rather reflects a curious fact of the word problem: the names in that list are all famous and well-known, but for reasons other than their research on the word problem! This is also true of G.\ S.\ Tseytin, whose contributions to theoretical computer science go far beyond his contributions to the word problem. Nevertheless, the beauty of his work on this problem, though limited to a single article, bears expanding on. The contribution comes to us in an essentially singular form: the construction of a semigroup with only seven defining relations and undecidable word problem. This is the semigroup $\fC$ with five generators $\{a, b, c, d, e\}$, and seven relations, which can easily be written on a blackboard:
\begin{align}
\nonumber ac &= ca, \quad ad=da, \quad bc=cb, \quad bd=db, \\ &eca = ce, \quad  edb = de, \quad cca = ccae.\label{Eq:Tseitin-semigroup} 
\end{align}
The fact that undecidability can be encoded into such a compact form may appear as a total mystery; it is the purpose of this article to both demystify this semigroup and demonstrate the beauty and importance of the ideas used to obtain this undecidability.

M.\ Dehn probably suspected that the word problem was decidable when he first posed it for groups in 1911 \cite{Dehn1911}, but he likely also suspected that it might be as hard as solving all of mathematics, and therefore for all practical purposes unsolvable; a stronger view was held by H.\ Tietze, who seems to have regarded problems of this form as genuinely unsolvable, see \cite[p.\ 563]{Stillwell2012}. Dehn died in 1952, the very same year that P.\ S.\ Novikov announced the algorithmic undecidability of the word problem for groups \cite{Novikov1952}; that year, one chapter of the history of group theory closed, another opened. Boone \cite{Boone1958} and Britton \cite{Britton1958}, both on the other side of the Iron Curtain, gave independent proofs of the same result a few years later. Eventually, Borisov \cite{Borisov1969} found an example of a 12-relator group with undecidable word problem. Occasionally, the history of the word problem for groups is summarized in a paragraph of approximately the same length and content as this. But this hides innumerable complexities and fascinating layers of the story, and a second reason for writing this article is also to unfold some less commonly discussed layers of this intricate story, while centering this on the semigroup constructed by Tseytin.

One such overlooked layer that will be particularly highlighted is the importance of semigroups in proving both the results of Novikov and Borisov. Semigroup theory is a subject which sometimes has a reputation as being rather dense and impenetrable, or at the very least as being filled with arcane notation and with results of little general interest (to outsiders). While it is indeed true that a cursory glance at Tseytin's article may give the impression that whatever he constructed is precisely of this nature, a somewhat longer glance reveals that nothing could be further from the truth. Instead, once the idea is laid bare, the proof is clear, and so too is the method by which it was arrived at. Furthermore, many of the techniques are today standard in combinatorics on words (such as factorization of words over a code), and many pages of technicalities and occasionally clumsy notation can be saved by eliminating the proofs of correctness of these techniques. We shall strive to do this as much as possible.

The article is structured as follows. First, we shall give the necessary background in \S\ref{Sec:background}, which will allow us to present Tseytin's semigroup $\fC$ and the proof of the undecidability of its word problem in \S\ref{Sec:Tseytin-proof}. Finally, in \S\ref{Sec:subsequents}, we will go through several ways in which echoes of Tseytin's proof and semigroup have appeared in the literature since then, including in the theory of special monoids, three-relation semigroups, the termination problem for rewriting systems, the word problem for groups, and universal semigroups. 

\section*{Acknowledgements}

The author wishes to thank Yu.\ V.\ Matiyasevich for encouragement to write this article, and for many useful notes, pointers to the literature, and comments leading to a significant improvement of the article.

\section{Background}\label{Sec:background}

Before presenting the contents of Tseytin's article \cite{Tseytin1958}, we will begin by laying its mathematical foundations in modern language; we will also, for the interested reader, occasionally mention the notation and terminology used by Tseytin. In particular, we shall assume no background of the reader with regards to semigroup theory. 

\subsection{Semigroup and monoid presentations} A \textit{semigroup} is a set endowed with a binary associative operation; if there is an identity element for this semigroup, then we call it a \textit{monoid}. Tseytin, and Soviet literature on the subject, for the most parts, calls monoids as \textit{associative calculi}. Let $A$ be any set, which we will call an \textit{alphabet}. We let $A^\ast$ denote the \textit{free monoid} on $A$, consisting of all words (i.e.\ finite-length sequences) over $A$, endowed with the operation of concatenation. This is a monoid, and any monoid is a homomorphic image of some free monoid. We let $A^+$ denote the \textit{free semigroup} on $A$, being the set of non-empty words over $A$ with concatenation. Letting $\varepsilon$ denote the empty word, the identity element of $A^\ast$, we obviously have $A^+ = A^\ast \setminus \{ \varepsilon \}$.

Every semigroup (resp.\ monoid) can be written as a quotient of some free semigroup (resp.\ monoid). Data encoding the kernel of the induced quotient map is called a \textit{presentation} for the monoid; more concretely, a semigroup presentation with generating set $A$ is a subset $R \subseteq A^+ \times A^+$ (resp.\ a subset $R \subseteq A^\ast \times A^\ast$), and is usually written 
\[
\pres{Sgp}{A}{R} \quad \text{(resp.\ $\pres{Mon}{A}{R}$).}
\]
The elements $(u, v) \in R$ are called \textit{defining relations} of $S$, and are usually written $u = v$ for simplicity. Any presentation gives rise to a semigroup (resp.\ monoid) in the following way; we only work with monoid presentations for simplicity. Let $\sim_R$ be the binary relation on $A^\ast \times A^\ast$ defined by: $x \sim_R y$ if and only if there is some relation $(u,v) \in R$ or $(v,u) \in R$ and words $x_0,x_1,y_0,y_1 \in A^\ast$ such that $x \equiv x_0ux_1$ and $y \equiv x_0vx_1$, i.e.\ if $x$ can be obtained from $y$ (or vice versa) by applying a single relation to $y$ (resp.\ $x$). Let $=_R$ be the reflexive, transitive closure of this relation. Then $=_R$ is an equivalence relation, and indeed a \textit{congruence} on $A^\ast$, i.e.\ the equivalence classes under $=_R$ form a monoid induced by multiplying representatives in $A^\ast$. We call this monoid $M = M(R)$. We shall say that two words $u, v \in A^\ast$ are \textit{equal} in $M$ if $u =_R v$. We usually write this as $u =_M v$, or say that ``$u=v$ in $M$''.

In this article, we shall be primarily concerned with monoids, rather than semigroups. There are two reasons for this: first, we shall require relations of the form $w=1$ in some of our presentations (which Tseytin calls \textit{special} monoids). Second, in any semigroup $S$ in which every defining relation is of the form $u=v$, where $u$ and $v$ are non-empty words, it follows easily that no word is equal to $1$. This is the case for Tseytin's semigroup (monoid) $\fC$ with undecidable word problem. Hence, the role of the identity element in such a case is purely as a formal symbol, and is occasionally useful in simplifying arguments, e.g.\ the base case of an induction. Hence, we frequently (intentionally) conflate semigroups and monoids throughout this present article.

We say that a monoid presentation $\pres{Mon}{A}{R}$ is \textit{finite} if $A$ and $R$ are both finite, and we say that a monoid $M$ is \textit{finitely presented} if it is isomorphic to the monoid of some finite presentation. We will exclusively deal with finitely presented monoids in this article. Henceforth, we shall blur the distinction between a monoid and a presentation for it, and write e.g.\ $M = \pres{Mon}{A}{R}$ to mean $M = M(R)$. Thus, for example, consider the monoid defined by the presentation 
\begin{equation}\label{Eq:monoid-ab-ba-pres}
M = \pres{Mon}{a,b}{ab=ba}.
\end{equation}
Here, the notation $ab=ba$ is used as shorthand indicate that $(ab,ba) \in R$ is a relation, and hence $ab \sim_R ba$. In this monoid, we have $ab =_M ba$, but also $a^m b^n =_M b^n a^m$ for all $m, n \in \N$, which follows by a direct induction. Indeed, this commutative monoid is isomorphic to $\N \times \N$, where the isomorphism is given by mapping any word $w$ to the pair $(\sigma_a(w), \sigma_b(w))$, where $\sigma_a$ (resp.\ $\sigma_b$) denotes the number of occurrences of the letter $a$ (resp.\ $b$) in $w$. 

Tseytin \cite{Tseytin1958} uses a somewhat different notation, but the same ideas: by $M : ab \bot ba$, it is meant that $ab=_M ba$ is an immediate consequence of one relation of $M$, i.e.\ $ab \sim_M ba$. For all $m, n \in \N$, Tseytin would have written $M : a^m b^n \dbot b^n a^m$ to say that $a^mb^n =_M b^n a^m$. However, we do \textbf{not} have $M : a^m b^n \bot b^n a^m$ for $m,n \neq 1$, as $a^m b^n = b^n a^m$ is not a defining relation of $M$.

\subsection{The word problem}

We will only define ``decision problem'' intuitively in this article, and use it to mean a problem of constructing an algorithm which, upon the input of some finite input data and a yes/no-question about this data, outputs the correct answer in a finite amount of time. If such an algorithm can be constructed, then we say that the decision problem is \textit{decidable}; otherwise, it is \textit{undecidable}. The reader interested in the precise formulations can consult e.g.\ \cite[Chapter~7]{Hopcroft1979}. 

Given a finitely presented monoid $M = \pres{Mon}{A}{R}$, there is a natural \textit{decision problem} associated to $M$ as follows. Let $u, v \in A^\ast$ be two words. Do we have $u =_M v$? Can we decide this algorithmically? This problem is called the \textit{word problem}\footnote{This problem has many names in the literature, a problem sometimes exacerbated by translations. It is often found, for example, as the \textit{identity, equality} or \textit{equivalence problem}; the latter is the translation I chose for the title of Tseytin's article \cite{Tseytin1958}.} for $M$. The word problem thus represents a fundamental question about the nature of $M$: is the most basic structure of $M$ -- whether or not we can tell apart elements in our chosen representation of them as words over $A^\ast$ -- computable from the finitely much information we have used to encode the (generally infinite) monoid $M$? For example, we certainly found that this was the case for $M = \pres{Mon}{a,b}{ab=ba}$ in \eqref{Eq:monoid-ab-ba-pres}. In other words, this monoid has decidable word problem. Do all finitely presented monoids have decidable word problem? This answer was first posed by A.\ Thue in 1914 \cite{Thue1914}, and was a question which was central to the entire development of algorithmic problems in algebra. We shall see in \S\ref{Subsec:decidability} that the answer to this question is a very firm \textbf{no}: there exist finitely presented monoids with undecidable word problem. Indeed, Tseytin proved that there exists a monoid with only seven defining relations and in which the word problem is undecidable. We shall present the idea of his proof in \S\ref{Sec:Tseytin-proof}. 

\subsection{Encodings} Let us first make a remark regarding the role of generators in a presentation versus that of relations. Any monoid with only one generator is cyclic, and there is little to say; however, once equipped with two generators, it is possible encode almost anything. For example, we may replace some finite number of generators $a_1, a_2, a_3, \dots$ by words over $\{ x, y\}^\ast$, for example by mapping $a_i$ to $xy^i$. Any word over $\{ xy^i \mid i \in \N \}$ can be uniquely decoded as a product of its generators, so we can recover the original word. Of course, this amounts to the statement that the free monoid on two generators contains embedded copies of all finitely generated (indeed countably generated) free monoids. This statement is familiar to anyone who has used a computer, the data of which is encoded as binary strings, i.e.\ elements of the free monoid $\{ 0, 1 \}^\ast$. In fact, by an analogous idea of encoding the generators appropriately over a two-generator alphabet, any finitely generated monoid with $k$ defining relations can also be embedded in a $2$-generated monoid with $k$ defining relations, as proved by Hall \cite{Hall1949}. This might seem like a prudent step that should always be performed to ``minimize'' a given presentation. However, this encoding process (from many generators to two) will tend to make the relations rather cumbersome and long, and the tension between having many generators and short relations needs to be resolved in order to find a truly minimal example. G. S. Tseytin's seven-relation semigroup is an example where this tension is minimized at five generators. 

One may also wish to encode the word problem of one monoid into that of another, even if the encoding is not effected by an embedding. For example, if one is given two monoids $M_i = \pres{Mon}{A_i}{R_i}$, where $i=1,2$, and one can say that 
\begin{equation}\label{Eq:M1-dec-M2-monoid}
u = v \text{ in $M_1$} \quad \iff \quad \varphi(u) = \varphi(v) \text{ in $M_2$},
\end{equation}
where $\varphi \colon A_1^\ast \to A_2^\ast$ is some computable function (not necessarily a homomorphism!), then this allows us to say that the word problem for $M_1$ can be reduced to the word problem for $M_2$. In other words, if we can solve the word problem in $M_2$, then we can solve the word problem in $M_1$; so if the latter has undecidable word problem, then so too does the former. Tseytin's idea was to show that there exists a seven-relation semigroup $\fC$ such that for every finitely presented group $G$, there exists a computable function $\varphi_G$ such that 
\begin{equation}\label{Eq:M1-dec-M2}
u = v \text{ in $G$} \quad \iff \quad \varphi_G(u) = \varphi_G(v) \text{ in $\fC$}.
\end{equation}
We emphasize the remarkable order of the quantifiers; the semigroup $\fC$ encodes \textit{every} finitely presented group. In particular, $\fC$ has undecidable word problem in the strongest possible sense. We remark briefly that Boone, Collins \& Matiyasevich \cite{Boone1974} in 1974 were able to construct a $\varphi_G$ of the above form with the additional property that it is an embedding, that the target semigroup has only four relations, and that $G$ may be taken as any finitely presented semigroup; we shall return to this point in \S\ref{Subsec:modern-proof-Matiyasevich}.

\subsection{Undecidability results}\label{Subsec:decidability}

The full history of the word problem for semigroups and groups, and the development of the same, is quite intricate and cannot possibly be told here in any degree of faithfulness. We will, however, give some indications. After the formalization of computability by Church, Turing, and many others, in the 1930s, it became natural to hunt for ``natural'' undecidable decision problems; all problems up to this point had been artificially created in the context of Turing machines, e.g.\ the printing problem (whether a machine will ever print a given symbol), or the halting problem (although, as pointed out by \cite{Lucas2021}, this problem was never studied by Turing himself). The word problem for finitely presented semigroups, studied by Thue in 1914, became a very natural target for an undecidable problem. 

In 1947, A.\ A.\ Markov \cite{Markov1947} and E.\ L.\ Post \cite{Post1947}, working independently of one another, proved the undecidability of the word problem in finitely presented semigroups; Markov's construction is more explicit, and from it one can extract a semigroup with 33 defining relations (and two generators) in which the word problem is undecidable. Many other natural problems were proved undecidable in this context in the next few years by Markov \cite{Markov1951}, including the isomorphism problem (given two monoid presentations, determine if the monoids presented by them are isomorphic). The next natural target was to try and prove the undecidability of the word problem for finitely presented \textit{groups}. But that undecidability would be the norm also for groups was far from a foregone conclusion. Indeed, B.\ H.\ Neumann announced, while working at the University of Manchester, likely in the early 1950s, that the word problem is decidable in every finitely presented group; another great mind present in Manchester at the time was Alan Turing -- who announced that the word problem is \textit{undecidable} in some finitely presented groups. Both later retracted their claims (!). But Turing was right, in the end, in two distinct ways: first, he published an article proving\footnote{Turing's article suffers from some errors, and an attempt to patch these was carried out by Boone \cite{Boone1958b}; the two articles should be read in parallel. Mindful of this, Novikov \& Adian \cite{NovikovAdian1958} later gave their own and somewhat different proof of the result.} the undecidability of the word problem in \textit{cancellative} semigroups (i.e.\ semigroups in which either of $xy=xz$ or $yx=zx$ imply $y=z$ for all $x,y,z$), and second, he was right: Novikov \cite{Novikov1952,Novikov1955} announced in 1952 the existence of groups with undecidable word problem, with a proof published in 1955 (although the group in the 1952 announcement differs from that in the 1955 article). Boone \cite{Boone1958} and Britton \cite{Britton1958} independently (of each other and of Novikov) gave proofs of the same theorem in 1958. A very readable and comprehensive treatment is given by Stillwell \cite{Stillwell1982}, which accomplishes the impressive task of both providing historical as well as mathematical details to an excellent degree. Rotman \cite[Chapter~12]{Rotman1995} also devotes a chapter to a proof that the word problem for groups is undecidable.

In the years following these developments, many new (too many to recount here) undecidability results appeared in the tailwaters of the word problem. Perhaps the most important in this context is the \textit{Adian--Rabin theorem} \cite{Adian1955, Rabin1958}, which was proved for monoids by Markov already in 1951 \cite{Markov1951}. Roughly speaking, this says that given a finite presentation of a group, we can say almost nothing about the properties of this group in general. For example, one cannot even decide if a finitely presented group is trivial or not! For more details, we refer the reader to the author's English translation \cite{NybergBroddaAdianRabin} of the articles on this subject by Adian and Markov. Notably, Tseytin \cite{Tseytin1956Markov} also contributed to this subject, providing a slight extension (focussed on the number of generators). The analogue of the Adian--Rabin theorem for cancellative semigroups was obtained in 1954 by Addison \cite{Addison1954} and, independently, Feeney \cite{Feeney1954}; somewhat more recently, it has also been extended to finitely presented \textit{inverse} semigroups by Yamamura \cite{Yamamura1997}. 

Thus, the early 1960s in group and semigroup theory was a world faced with a sudden abundance of undecidability. It was in this world of newfangled impossibility that Tseytin discovered his seven-relation semigroup with undecidable word problem.

\section{Tseytin's proof}\label{Sec:Tseytin-proof}

In this section, we give an overview, in modern terminology, of how Tseytin's proof works. Tseytin's proof and article are both complete and essentially self-contained, and modern proofs of stronger results are available elsewhere (see \S\ref{Subsec:modern-proof-Matiyasevich}), so we shall only attempt to convey the high-level idea here. All references to sections in Tseytin's article \cite{Tseytin1958} will, for clarity, be accompanied also by the citation, e.g.\ as \cite[\S2]{Tseytin1958} rather than simply \S2. In other words, all section numbering without such a citation refers to a section within this present article. We remark that Tseytin \cite{Tseytin1956}\footnote{The author wishes to thank L. D. Beklemishev and A. L. Talambutsa for providing him with a copy of this article.} announced his result already in 1956, although the semigroup therein has slightly different relations, and no proof is provided (the outline given shows that it doubtlessly uses the same ideas). We therefore concern ourselves only with the 1958 article. As mentioned in the introduction, an English translation of \cite{Tseytin1958} by the author of this present article is supplemented at the end of this present article. 

\subsection{Codes}

The first mathematical section \cite[\S2]{Tseytin1958} of Tseytin's article contains only elementary results in combinatorics on words. Let $A$ be an arbitrary finite alphabet. Let $\sigma : A \to \N$ be any injective function. The main tool introduced is an encoding function\footnote{This function is essentially the inverse function to the ``compression'' function introduced by Adian \& Oganesian \cite{Adian1978} in their study of the word problem for one-relation monoids; see \cite{Zhang1992Compression, NybergBroddaCompression} for a modern treatment of this via rewriting systems.}, which is a function $\varphi_\circ : A^\ast \to \{ a, b\}^\ast$ defined as follows. We have $\varphi_\circ(\varepsilon) = a$, and if $w = a_1 a_2 \cdots a_n$ where $a_i \in A$, then 
\begin{equation}\label{Eq:def-phi-circ}
\varphi_\circ(w) = ab^{\sigma(a_1)} ab^{\sigma(a_2)} \cdots ab^{\sigma(a_n)}a.
\end{equation}
Note that $\varphi_\circ$ is not a homomorphism. However, it is clear that $\varphi_\circ$ is injective on non-empty words; given a word written of the form of the right-hand side of \eqref{Eq:def-phi-circ}, we can find what word $w$ gave rise to it, the only possible ambiguity being taking $w = a_i$ where $a_i \in A$ is a letter such that $\sigma(a_i) = a$, for then $\varphi_\circ(a_i) = a = \varphi_\circ(\varepsilon)$. This uniqueness statement amounts, effectively, to saying that the set $I = \{ ab^i \mid i \in \N \}$ freely generates a free submonoid of $\{a, b\}^\ast$, i.e.\ $I$ is a code. Proving this is the purpose of the first few lemmas of \S2.

Let $S = \pres{Mon}{A}{u_i = v_i}$ be any finitely presented monoid (we do not index the relations for brevity). Then we can define the monoid $S^\circ$ by 
\[
S^\circ = \pres{Mon}{a,b}{\varphi_\circ(u_i) = \varphi_\circ(v_i)}.
\]
Then (Lemma~\S2.6 and Lemma~\S2.9) it is not hard to show that for any two words $w_1, w_2 \in A^\ast$, we have $w_1 =_S w_2$ if and only if $\varphi_\circ(w_1) =_{S^\circ} \varphi_\circ(w_2)$. That is, $\varphi_\circ$ can be used to encode the monoid $S$, which may be generated by many elements, into a $2$-generated monoid $S^\circ$. This is an easy consequence of the code property of $I$ above.\footnote{The interested reader may compare this with \cite[Theorem~3]{Adian1978}, especially equation (8) of the same.} 

\subsection{Encoding groups}

In \S3, the semigroup $\fC$ with undecidable word problem is introduced. This depends on a natural number parameter $i \geq 0$, which will correspond to a choice for an encoding function $\varphi_\circ$. The semigroup $\fC^{(i)}$ is then defined as the semigroup with five generators $\{ a,b , c, d, e\}$ and seven relations
\begin{align}
\nonumber ac = ca&, \quad ad=da, \\ \nonumber bc=cb&, \quad bd=db, \\ eca &= ce, \tag{\ref{Eq:Tseitin-semigroup}} \\ \nonumber edb &= de, \\ \nonumber cd^i ca &= cd^i cae.
\end{align}
Throughout this present discussion (and, as the observant reader will note, in the introduction to this article), we shall take $i=0$ for ease of notation, and write $\fC = \fC^{(0)}$. The proof carries through for all $i$. For our discussion of the specific word problem in \S\ref{Subsec:Specific-word-problem}, in order to minimize the presentation we will need to take $i=1$. 

For our exposition, we will now change the names of the generators of $\fC = \fC^{(0)}$. We shall also refer to it as a monoid. The monoid $\fC$ has $5$ generators, which we denote by $\{ x_1, x_2, y_1, y_2, t \}$. Let $F_x$ be the free monoid on the two generators $\{x_1, x_2\}$, and $F_y$ the free monoid on $\{y_1, y_2 \}$. Take the direct product\footnote{The direct product $F_2 \times F_2$ of two free monoids or groups is a frequent culprit in the construction of undecidable problems in algebra; for example, the direct product of two free groups is proved by Mikhailova \cite{Mikhailova1958} to have undecidable membership problem, which leads to analogous results in $\operatorname{SL}_4(\Z)$. Similarly, the direct product of two free monoids is frequently used in recent decision problems for matrix groups, cf.\ e.g. \cite{Ko2018}.} $F_x \times F_y$ and add three relations: 
\begin{align*}
ty_k x_k &= y_k t, \quad (k = 1, 2)  \\
y_1^2 x_1 &= y_1^2 x_1 t.
\end{align*}
We call the resulting seven-relation monoid $\fC$, i.e.\
\begin{equation}\label{Eq:fC-def-txy}
\fC = \pres{Mon}{x_1,x_2,y_1,y_2, t}{x_k y_\ell = y_\ell x_k, \: t y_k x_k = y_k t, \: y_1^2 x_1 = y_1^2 x_1 t},
\end{equation}
where the indices on the relations range over all $1 \leq k, \ell \leq 2$. Of course, this is the same monoid as that defined by \eqref{Eq:Tseitin-semigroup} with $i=0$.  

Let us give some intuition for this monoid \eqref{Eq:fC-def-txy}. A useful analogy is to think of $\fC$ as a record player: when we give it a record, which encodes the defining relations of an arbitrary group, it can simulate the insertion and deletions of the relations of this group without altering the record. The letter $t$ will act like the ``stylus'' of the record player, running back and forth to insert and delete relations in appropriate places. This allows $\fC$ to simulate any group within itself. Slightly more formally, let 
\[
G = \pres{Gp}{x_1,x_2}{R_1 = 1, \dots, R_k = 1}
\]
be a group, where the notation $\pres{Gp}{A}{R}$ is used to indicate that we are dealing with a \textit{group} presentation, i.e.\ for all $a \in A$, we also have generators $a^{-1}$, and the relations $aa^{-1} = a^{-1} a = 1$. We will assume some strong properties about the presentation: we require that all $R_j$ are words containing no inverse symbols, and that $G$ coincides with the monoid defined by the same presentation. This can always easily be done. From such a presentation of $G$, we will construct a record $\bS = \bS(G)$, which should be though of as the word $R_1 \mid R_2 \mid \cdots \mid R_k$ of the relators of the group separated by some separating symbols (themselves encoded). The words $R_j$ in $\bS$ will also be written over $y_1, y_2$ rather than $x_1, x_2$, by swapping every $x_j$ with $y_j$ for $j=1,2$. Then, given any two words $u, v \in F_x$, the main technical result \cite[Lemma~\S4.5, Lemma~\S5.12]{Tseytin1958} states:
\begin{equation}\label{Eq:main-lemma}
u = v \text{ in $G$} \quad \iff \quad \bS u = \bS v \text{ in $\fC$}.
\end{equation}
This, once properly stated, will give the desired result: it will follow that the word problem in $\fC$ is decidable only if the word problem is decidable in \textit{every} finitely presented group $G$. That is, $\fC$ can be used to simulate every finitely presented group. But the word problem is far from being decidable in every finitely presented group: indeed, there exist finitely presented groups which have undecidable word problem! This shows that the word problem in $\fC$ is (very) undecidable. 

\subsection{Details of the encoding}\label{Subsec:details-proof}

We now add even further details to the above outline. First, it is not clear how to encode a given group $G$ into the word $\bS$ in a clear manner. Let us describe this process in some more detail. 

First, in order not to deal with inverse symbols, we will instead work with \textit{special monoids} rather than groups. A monoid $M$ is said to be \textit{special} if the right-hand side of every relator is always $1$, i.e.\ if every relation is of the form $R = 1$. Obviously, every group is special, but not every special monoid is a group; the canonical non-example is the \textit{bicyclic} monoid $\pres{Mon}{b,c}{bc=1}$, in which there are no non-trivial invertible elements ($b$, for example, is right invertible, but not left invertible). 

We will encode a given special monoid into a single word $\bS$ as follows. Let 
\[
M = \pres{Mon}{A}{R_1 = 1, R_2 = 1, \dots, R_k = 1}
\]
be a finitely presented special monoid. Let $x$ be a new symbol. We let $M_1$ be the special monoid with the generating set $B = A \cup \{ x \}$ and the relators $R_j \alpha = 1$ (for all $j$), together with a relator $\alpha = 1$. Of course, $M_1 \cong M$, indeed $M$ embeds into $M_1$ by the identity mapping. This is the content of \cite[\S3, Lemmas~1--4]{Tseytin1958}. Taking an encoding function 
\[
\varphi_\circ \colon B^\ast \to F_x
\]
defined as in \eqref{Eq:def-phi-circ} but with the constraint that $\varphi_\circ(\varepsilon) = x_1$ and $\varphi_\circ(x) = x_1 = x_1x_2^0$. (We note that it is here that the choice of $i \in \N$ in $\fC^{(i)}$, which we took as $i=0$ for simplicity, comes into play, by choosing $\varphi_\circ(x) = x_1x_2^i$.) Then $M_1^\circ$ is well-defined; in particular all its relations are of the form $\varphi_\circ(R_i) = x_1$. 

\begin{example}\label{Ex:bicyclic}
If $M = \pres{Mon}{b,c}{bc=1}$, the bicyclic monoid, then 
\begin{align*}
M_1 &= \pres{Mon}{b,c,x}{xbxcx = 1, \: x = 1}, \\ 
M_1^\circ &= \pres{Mon}{x_1,x_2}{x_1 (x_1x_2) x_1 (x_1x_2^2) x_1 x_1 = x_1, \: x_1^2 = x_1},
\end{align*}
taking $\varphi_\circ(b) = x_1x_2, \varphi_\circ(c) = x_1x_2^2$. Of course, $M_1 \cong M$. 
\end{example}

Continuing the proof, one must show \cite[\S3, Lemmas 5 \& 7]{Tseytin1958} that $\varphi_\circ(u) = \varphi_\circ(v)$ in $M_1^\circ$ if and only if $u = v$ in $M$. This is the final result of \cite[\S3]{Tseytin1958}. Thus $M$ and $M_1^\circ$ are essentially equivalent; we have encoded our special monoid $M$ into a suitable monoid $M_1^\circ$, which is no longer special.\footnote{In the terminology of Kobayashi \cite{Kobayashi2000}, $M_1^\circ$ is \textit{subspecial}, and one checks that $M_1^\circ$ \textit{compresses} to its \textit{left monoid} $M_1$. Indeed, it follows from the theory of compression, e.g.\ in Lallement \cite{Lallement1974} or Adian \& Oganesian \cite{Adian1978}, that the word problem for $M_1^\circ$ is equivalent to that of $M_1 \cong M$.}

Next, the crucial step: it is shown that $\fC$ can simulate the word problem of $M_1^\circ$, in the sense of \eqref{Eq:M1-dec-M2}. At the end of \cite[\S3]{Tseytin1958}, a word $\bS$ is introduced, which encodes the defining relations of $M_1^\circ$; it has the form 
\[
\bS = \sigma [\varphi_\circ(x R_1 x R_2 x \cdots R_k x)].
\]
Here $\sigma$ is the function which swaps $x_i$ for $y_i$, i.e.\ the morphism $\sigma \colon F_x \to F_y$ extending $x_1 \mapsto y_1$ and $x_2 \mapsto y_2$. In particular, $\bS \in F_y$. This word $\bS$ now contains all of the information about the presentation of $M_1^\circ$, which in turn is a simple encoding of $M$. As we shall now see, the relations of $\fC$ will be able to use this information to simulate the word problem of $M_1^\circ$. 

\begin{example}
Continuing Example~\ref{Ex:bicyclic}, we must encode the relations of $M_1^\circ$ into a single word: 
\begin{align*}
\bS = \sigma[\varphi_\circ(x R_1 x R_2 x \cdots R_k x)] &= \sigma[\varphi_\circ(xbcx)] \\
&= \sigma[x_1 (x_1x_2) x_1 (x_1x_2^2) x_1 x_1] \\
&= y_1(y_1y_2)y_1(y_1y_2^2)y_1y_1.
\end{align*}
Thus, the word $y_1^2y_2y_1^2y_2^2y_1^2$ is a record; when it is given to the record player $\fC$, the music that is the word problem of the bicyclic monoid will begin to play.
\end{example}

Continuing with the proof, what needs to be proved is the following crucial statement: if $u, v \in \{ x_1, x_2 \}^\ast$, then:
\begin{equation}\label{Eq:main-claim}
u = v \text{ in } M_1^\circ \quad \iff \quad \bS u = \bS v \text{ in } \fC \tag{$\ast$}
\end{equation}
It is clear that from this statement, everything which was needed for undecidability follows; indeed, as any group $G$ is a special monoid, we will have established \eqref{Eq:main-lemma}\footnote{In \S\ref{Subsec:specialmonoids} of this present article, we shall see that \eqref{Eq:main-claim} is not, in fact, much stronger than \eqref{Eq:main-lemma}, even though not every special monoid is a group.}. The direction $(\implies)$ of \eqref{Eq:main-claim} is proved in \cite[\S4]{Tseytin1958}, specifically in Lemma~5 of the same. The direction $(\impliedby)$ is proved in \cite[\S5]{Tseytin1958}, specifically Lemma~12 of the same.

The proof of $(\implies)$ is beautiful. Let us show how to use $\fC$ to apply the defining relation $(\varphi_\circ(R_k x) = x_1)$ of $M_1^\circ$ to the word $u$. Specifically, we will replace an $x_1$ inside $u$ by $\varphi_\circ(R_k x)$. It is here that the ``stylus'' nature of $t$ comes into play. Now, we may assume $u$ contains an occurrence of $x_1$. This is because every relation of $M_1^\circ$ contains an $x_1$. Write $u \equiv u_0 x_1 u_1$, where $u_0, u_1 \in F_x$. Since $\bS$ ends with $\sigma(\varphi_\circ(xR_k x)) \in F_y$, and since this word commutes with every element of $F_x$, it follows that we can, as an equality of words, write
\[
\bS u = \bS u_0 x_1 u_1 = \bS' u_0 \sigma(\varphi_\circ(x R_k x)) x_1 u_1
\]
where $\bS = \bS'\sigma(\varphi_\circ(xR_kx))$. Now $\sigma (\varphi_\circ(xR_k x))$ begins and ends with $y_1^2$. Hence, using the relation $y_1^2x_1 = y_1^2 x_1 t$, we have that 
\[
\bS' u_0 \sigma(\varphi_\circ(xR_k x)) x_1 u_1 = \bS' u_0 \sigma(\varphi_\circ(R_k x)) x_1 t u_1 \quad \text{ in $\fC$.}
\]
Let $\varphi_\circ(x R_k x) = x_1^2 R$, where $R \in F_x$. Then 
\[
\sigma(\varphi_\circ(x R_k x))x_1t = \sigma(x_1^2 R) x_1 t = y_1^2 x_1  \sigma(R) t \quad \text{ in $\fC$,}
\]
as $x_1$ commutes with every element of $F_y$. The next step is crucial. The relations $t y_j x_j = y_j t$ and the relations of $F_x \times F_y$ now clearly imply (\cite[\S4, Lemma~3]{Tseytin1958}) that $\sigma(X) t = t \sigma(X) X$ for any word $X \in F_x$, with the letter $t$ moving as a stylus over $\sigma(X)$, reading it from right to left, and writing down a copy of it, but written over the alphabet $F_x$. Thus, applying this process, moving the stylus $t$ across the word $\sigma(R)$, we find
\[
y_1^2 x_1 \sigma(R) t = y_1^2 x_1 t \sigma(R) R = y_1^2 x_1 \sigma(R) R = y_1^2 \sigma(R) x_1 R \quad \text{ in $\fC$.}
\]
The second equality uses the relation $y_1^2x_1 t = y_1^2x_1$, and the last uses the relations of $F_x \times F_y$ to permute $x_1 \in F_x$ with $\sigma(R) \in F_y$. By the properties of the encoding function, it is not hard to verify that (as words) we have
\[
y_1^2 \sigma(R) x_1 R = \sigma(\varphi_\circ(x R_k x)) \varphi_\circ(R_k x)
\]
and hence we have completed our insertion of $\varphi_\circ(R_k x)$. We now retrace our steps by returning the word $\sigma(\varphi_\circ(R_k x))$ to the end of $\bS'$ whence we took it. In this way, we obtain the word $\bS u_0 \varphi_\circ(R_k x) u_1$, i.e.\ 
\[
\bS u_0 x_1 u_1 = \bS u_0 \varphi_\circ(R_k x) u_1 \quad \text{in $\fC$.}
\]
Thus we have applied the relation $x_1 = \varphi_\circ(R_k x)$ of $M_1^\circ$ to the word $\bS u = \bS u_0 x_1 u_1$. Of course, this could also be done in reverse; and, without writing out the technical details, this can be checked to work for any $R_i$, not just $R_k$. In other words, we can use $\fC$ to apply the relations of $M_1^\circ$ to the word $u$ inside $\bS u$. This gives us the direction $(\implies)$. Some parts of this proof could fruitfully, in a more detailed proof, be interpreted as the letter $t$ acting on pairs of words in $F_x \times F_y$; for example, we have used the fact (\cite[\S4, Lemma~3]{Tseytin1958}) that $\sigma(u)t = t \sigma(u)u$ in $\fC$ for all $u \in F_x$, which can be interpreted as an action
\[
(1, \sigma(u)) \cdot t = t \cdot (u, \sigma(u)) 
\]
of $t$ on $F_x \times F_y$. We leave the interested reader to develop this further.

For the direction $(\impliedby)$, covered in \cite[\S5]{Tseytin1958}, the key idea is to ensure that no other types of word transformations can be applied to the $\bS u$ other than precisely those which arise from $M_1^\circ$. For this, we will make one observation: as noted above, the word $\bS$ will be a word in $F_y$. The relations of $F_x \times F_y$ can hence not change anything about $\bS$ (other than permuting some of its letters with letters from $F_x$, which can always be undone). Similarly, the relations $ty_j x_j = y_j t$ will not modify $\bS$; neither will the relation $y_1^2 x_1 = y_1^2 x_1 t$. Consequently, the ``record part'' $\bS$ of the word $\bS w_1$ will remain unchanged when we have a word $w_1 \in F_x$, and any word equal to $\bS w_1$ will contain a copy of $\bS$. This is one part of the argument, and the remaining part follows from a precise analysis of where relations of $\fC$ can be applied to words of the form $\bS u$. This is an exercise in bookkeeping, and requires no real new ideas (unlike the proof of the forward implication). We therefore omit this bookkeeping; and thus, by the miracle of omission, we have proved both directions of \eqref{Eq:main-claim}. This completes the proof of the undecidability of the word problem in $\fC = \fC^{(0)}$. More generally, it also shows undecidability of the word problem in $\fC^{(i)}$ for all $i \geq 0$.

\subsection{Specific word problem}\label{Subsec:Specific-word-problem}

The word problem is undecidable in finitely presented groups; in particular, there exists a finitely presented group $G$ such that there is no algorithm which takes as input a word $w$ in the generators of $G$ and decides whether or not $w = 1$ in $G$. This is a ``one-sided'' word problem, in the sense that it only takes one word as an input, and compares it to a fixed word ($1$, in this case). In a group, this is clearly equivalent to the ``two-sided'' word problem, which takes as input \textit{two} words $u, v$ and decides whether or not $u = v$ in $G$ (for $u = v$ if and only if $uv^{-1} = 1$). But for monoids in general, this trick does not work. Nevertheless, $\fC$ has undecidable ``one-sided'' word problem, too, which we will now show. 

We begin with some proper definitions. Let $M = \pres{Mon}{A}{R}$. The \textit{specific} word problem for $M$ with respect to a fixed word $w \in A^\ast$ is the problem of deciding, on input $u \in A^\ast$, whether or not $u = w$ in $M$. The specific word problem is clearly equivalent to the word problem in a group, but this is not necessarily so in a monoid; there may be a solution to the specific word problem for any fixed word, but that this solution is not uniform. However, Tseytin's semigroup $\fC$ has undecidable specific word problem, too. Namely, fix a group $M$, written as a special monoid, with undecidable word problem. If one now takes the fixed word $\bS x$, where $\bS$ encodes the relations of $M$, then solving the specific word problem with respect to this word is, by \eqref{Eq:main-lemma}, the same as asking whether a word is equal to $x$ in $M_1^\circ$, which is equivalent to whether or not a word is equal to $1$ in $M$, which is undecidable by assumption. However, $\bS x$ is, in general, a very long word, since even the smallest known presentation for an $M$ with undecidable word problem has many relations.

Tseytin rectifies this by constructing another semigroup, which he calls $\fC'$, in \cite[\S7]{Tseytin1958}. This depends on two natural number parameters $i, j \in \N$, with the crucial property that $i \neq j$. We will use the notation $\fC^{(i,j)}$ to be the semigroup with the same generators and relations as $\fC^{(i)}$ (from \S\ref{Subsec:details-proof}), together with two more relations
\begin{equation}\label{Eq:new-relations-general}
cab^jab^ja = ab^jab^ja, \quad dab^jab^ja = ab^jab^ja.
\end{equation}
Then Tseytin proves the following remarkable theorem: in the nine-relation semigroup $\fC^{(i,j)}$, there is no algorithm for deciding whether or not a given word is equal to $ab^jab^ja$, i.e.\ the specific word problem is undecidable with respect to this word. 

To minimize the total length of the relations and the specific word, we will choose $i=1$ and $j=0$ (in \S\ref{Subsec:details-proof}, we chose $i=0$). Thus, we have a nine-relation semigroup $\fC^{(1,0)}$ with the relations of $\fC^{(1)}$ together with 
\begin{equation}\label{Eq:new-relations}
ca^3 = a^3, \quad da^3 = a^3.
\end{equation}
Throughout this section, we set $\fC' = \fC^{(1,0)}$ (Tseytin uses the notation $\fC_2$ for this particular semigroup). In particular, for this choice of $\fC'$, and in the notation of \eqref{Eq:fC-def-txy}, these new relations become $y_1x_1^3 = x_1^3, \: y_2x_1^3 = x_1^3$. We will now sketch a proof that the specific word problem with respect to the word $a^3$ in $\fC'$ is undecidable.

We only give a brief mention of the idea behind this proof. One fixes a special monoid $M = \pres{Mon}{A}{R}$, as before, and defines analogous encodings. The idea is straightforward: using the relations \eqref{Eq:new-relations}, written as $y_1x_1^3 = x_1^3, \: y_2x_1^3 = x_1^3$, one can ``wipe away'' the word $\bS$ from $\bS x$, once one is done using it. This is accomplished by adding a new letter $y$ (denoted $\beta$ in \cite{Tseytin1958}). We fix the encoding $\varphi_\circ$ as before, and such that $\varphi_\circ(y) = x_1$. It is at this point that our choice of $i=1$ in the relation $cd^ica = cd^icae$, and more generally that $i \neq j$, in order to ensure that the encoding does not break; indeed, we will have $\varphi_\circ(x) = x_1x_2^i = x_1 x_2$ for our choice of $i$. We consider the free product $M' = M \ast \langle y \rangle \cong M \ast \mathbb{N}$. After some elementary lemmas about free products, one shows the following main claim. Let $P$ be any word over $A$ (the alphabet of $M$). Then we have:
\begin{equation}\label{Eq:main-claim-particular}
P = 1 \text{ in } M' \quad \iff \quad \bS \varphi_\circ(y P y) = x_1^3 \text{ in } \fC'. \tag{$\ast \ast$}
\end{equation}
The forward direction ($\implies$) is rather simple to prove, and is done in \cite[\S7, Lemma~7]{Tseytin1958}. By analogous reasoning to the previous sections, and using the relations of $\fC^{(i)}$ inside $\fC'$, one shows that $\bS \varphi_\circ(y P y) = \bS \varphi_\circ(yy)$ in $\fC'$. Since $\varphi_\circ(yy) = x_1^3$, we can now apply the new relations $y_1x_1^3 = x_1^3, \: y_2x_1^3 = x_1^3$ to the word $\bS \varphi_\circ(yy) = \bS x_1^3$, and since $\bS \in F_y$ it follows that $\bS x_1^3 = x_1^3$ in $\fC'$. This proves the forward direction. The converse direction, proved in \cite[\S7, Lemma~8]{Tseytin1958}, is not much harder to prove, and relies on the control obtained from the previous section together with the fact that no word equal to $1$ in $M'$ can contain $y$ as a subword, which follows from the free product decomposition $M' \cong M \ast \langle y \rangle$. Having proved the main claim \eqref{Eq:main-claim-particular}, it is now clear that the result follows: since no algorithm can decide whether a given word is equal to $1$ in $M \leq M'$, it follows from \eqref{Eq:main-claim-particular} that no algorithm can decide whether or not a given word is equal to $x_1^3$ in $\fC' = \fC^{(1,0)}$. More generally, it also shows undecidability of the specific word problem with respect to $x_1x_2^jx_1x_2^jx_1$ in $\fC^{(i,j)}$ for all non-negative $i \neq j$.

\subsection{D.\ Scott's semigroup}

We remark that Dana Scott \cite{Scott1956} also produced a seven-relation semigroup with undecidable word problem, and this was announced already in 1956, just like Tseytin's \cite{Tseytin1956}. No proof was published outside of the cited announcement in the \textit{Journal of Symbolic Logic}, but the semigroup bears a strong resemblance to $\fC$: there are six generators $a,b,c,d,e,f$, and ten relations:
\begin{align*}
ac = ca, \: ad &= da, \: ae = ea, \\
bc = cb, \: bd &= db, \: be=eb, \\
e = ef, \: e = fe, &\: fca = cf, \: fdb = df
\end{align*}
which, then, by a ``suitable encoding'' can be reduced to seven; whereof four are encodings of the first six, one is an encoding of $e = ef = fe$, and the last two remain unchanged. Undoubtedly, the idea of the proof is similar, though no proof was ever published. D.\ Scott (personal communication) has informed me that he became aware of the fact that the same result had appeared in the USSR, and that this was the primary reason for foregoing publication. It is an interesting exhibition of mathematical symmetry that the knowledge of Scott's semigroup also made it across the Iron Curtain, albeit in the other direction: it is referenced, for example, by Makanin \cite{Makanin1966} in an article in the \textit{Doklady Akademii Nauk} in 1966.

\section{Subsequent developments}\label{Sec:subsequents}

In this section, we will describe some ways in which the article and ideas introduced by G.\ S.\ Tseytin influenced research on combinatorial group and semigroup theory in the years and decades following their appearance. 

\subsection{Special monoids}\label{Subsec:specialmonoids}

As we have seen in \S\ref{Sec:Tseytin-proof}, one of the key constructions used in encoding groups with undecidable word problem into $\fC$ was that of \textit{special monoids}, i.e.\ monoids with a presentation in which all defining relations are of the form $r = 1$. These were first defined by Tseytin, but would be studied in much greater depth a few years later by S.\ I.\ Adian in his seminal 1966 monograph \cite{Adian1966}. For example, Adian proved that a special monoid is finite if and only if it is a group. Perhaps the most important result, however, related to the group of units (i.e.\ the subgroup consisting of all two-sided invertible elements) and the word problem. This story was completed by G.\ S.\ Makanin \cite{Makanin1966}, with proofs given in his 1966 Ph.\ D. thesis \cite{Makanin1966Thesis} (translated by the author in \cite{NBMakanin1966Thesis}). He extended the results by Adian to all finitely presented special monoids, yielding the following theorem:

\begin{theorem*}[Adian \& Makanin, 1966]
Let $M$ be a $k$-relation special monoid, and let $G$ be its group of units. Then:
\begin{enumerate}
\item $G$ is a $k$-relator group;
\item The word problem for $M$ reduces to that of $G$. 
\end{enumerate}
\end{theorem*}

A modern proof, using the language of string rewriting, has been given by Zhang \cite{Zhang1992}. The author of the present article has also studied special monoids from the point of view of formal language theory in \cite{NybergBroddaSpecial}. 

In light of the above theorem, one note, given at the end of \S6 of Tseytin's article, stands out. Paraphrased, this states that if the (uniform) word problem for special monoids can be proved undecidable in a direct manner, then the proof of the undecidability of the word problem in $\fC$ can be greatly simplified (as it would bypass the need to use Novikov's group with undecidable word problem). However, Tseytin then states that he is not aware of such a direct proof. Indeed, the theorem above by Adian \& Makanin directly shows that no such direct proof is feasible: given any $k$-relation special monoid with undecidable word problem, there exists a $k$-relator group with undecidable word problem (the converse is not necessarily true; see the next paragraph). Thus, Tseytin's introduction of special monoids in this context is commendable, as it introduced no additional strong overhead on the complexity, while simultaneously simplifying the proof. 

We make a final remark about special monoids and how they relate to combinatorial group theory. Of course, every $k$-relation special monoid which is a group is defined by the group presentation with the same generators and relations, and is therefore necessarily also a $k$-relator group. However, not every $k$-relator group is a $k$-relation special monoid. A classical example of this is $\mathbb{Z}^2$, which is isomorphic to the one-relator group $\pres{Gp}{a,b}{aba^{-1}b^{-1}=1}$, but which admits no one-relation monoid presentation (special or otherwise); a simple proof of this fact is given by Otto \cite{Otto1988}. A one-relator group which admits a one-relation monoid presentation is called \textit{positive}; it is (non-trivially) equivalent to the group admitting a presentation of the form $\pres{Gp}{A}{w=1}$, where $w$ is a positive word (see Perrin \& Schupp \cite{Perrin1984}). Baumslag \cite{Baumslag1971} proved that all positive one-relator groups are residually solvable (i.e.\ the intersection of all terms in the derived series is trivial). Recent results by the author, Foniqi \& Gray \cite{IsoGrayCF}, however, have demonstrated that even positive one-relator groups can exhibit rather pathological behaviour. The above result by Adian, proved in the case of a single relation in \cite{Adian1966}, is actually somewhat stronger: it proves that the group of units of a one-relation special monoid is a \textit{positive} one-relator group, and furthermore gives a simple procedure for computing it. We refer the reader to \cite[\S2.3]{NybergBroddaSurvey} for a detailed account.

\subsection{Five relations}\label{Subsec:five-relations}

Few encodings remain optimal for long. This was also the case of the semigroup constructed by Tseytin, and within 10 years improvements and encodings of Tseytin's seven-relation semigroup $\fC$ were discovered. 

The first improvement was rather modest, but turned out to be useful. This was an observation by Yu. V. Matiyasevich \cite{Matiyasevich1967b} that the final relation $cca=ccae$ in $\fC$ can be replaced with $cca = cce$ with no loss of undecidability, and that Tseytin's proof carries through (in the context of the \textit{specific} word problem, this was discussed already in \S\ref{Subsec:Specific-word-problem}). The resulting semigroup has a total of $31$ occurrences of letters in its defining relations (compared to Tseytin's $32$) and remains to this day the smallest, with respect to this particular metric, semigroup with undecidable word problem. 

However, the main contribution in the article \cite{Matiyasevich1967b} is not this reduction, but something more significant: the construction of a five-relation semigroup with undecidable word problem. This uses the above improvement, and overall the construction is quite analogous to Tseytin's, by further encoding the obtained relations. The semigroup (see \cite[p. 56, Theorem~2]{Matiyasevich1967b}) has two generators $x,y$, and the following five defining relations:
\begin{align*}
xyx^2y^2 &= y^2x^2yx, \quad x^2yxy^2x = y^2x^3yx, \\
xyx^3y^2 &= xy^2xyx^2, \:\: x^4y^2x^2yx = y^2x^4, \\
& y^3x^2y^2x^2yx = y^3x^2y^2x^4. 
\end{align*}
The proof relies on Tseytin's $\fC$, and uses a clever encoding trick to compress the seven relations into five. We do not go into the details of this encoding here, due to the three-relation construction we will present in \S\ref{Subsec:three-relations}. 

Around the same time, G. S. Makanin announced, in an addition to the \textit{Doklady} article \cite{Makanin1966}, that he had proved the same result (while acknowledging that Matiyasevich had already proved it); Makanin's semigroup with undecidable word problem, added in proof to \cite{Makanin1966}, consists of three generators $x,y,z$ and the following five relations:
\begin{align*}
z^2 y^2 &= y^2z^2, \quad yz^3y^2 = z y^3 z^2, \\
xz^2y^2 &= y^2x, \quad\: xy z^3 y^2 = z y^2 x, \\
& y^2 z^2 y^4 z^2 = y^2 z^2 y^4 z^2 x.
\end{align*}
Unlike Matiyasevich's example, no proof was ever provided of the undecidability of the word problem in this semigroup, nor even an indication of the method by which it was obtained. Yu.\ Matiyasevich (private communication) has informed the author of what encoding was used, and it is included here as it does not appear anywhere else in the literature: first, in Tseytin's semigroup $\fC = \fC^{(0)}$ one replaces the relation $cca = ccae$ by $cccaa = cccaae$, which does not change the undecidability result. Next, one encodes by $a \mapsto zz, b \mapsto yzz, c \mapsto yy, d \mapsto zyy$, and $e \mapsto x$. Using this mapping, it is not hard to see that two of the relations become redundant, yielding the above presentation.

\subsection{Three and fewer relations}\label{Subsec:three-relations}

The step from seven to five relations represents a modification of the idea by G. S. Tseytin, and is a local improvement of encoding ideas. By contrast, the next step, also carried out by Matiyasevich, is very innovative. This was his 1967 construction of a three-relation semigroup with undecidable word problem. 

For the other semigroups with undecidable word problem in this article (except those constructed by Markov and Post), we have written down explicit presentations. For Matiyasevich's semigroup, denoted $\mathfrak{N}$ in \cite{Matiyasevich1967}, this will not quite be the case. The presentation has two generators $a, b$, and the following three relations:
\begin{equation}\label{Eq:matiyasevich-semigroup}
aabab = baa, \quad aabb = baa, \quad W_1 = W_2
\end{equation}
The first two are easy enough to commit to memory. The third relation $W_1 = W_2$ is not: the word $W_1$ contains $304$ letters, and the word $W_2$ contains $608$ letters. Thus, not only is it necessary to encode the word problem for all groups into the relations of a semigroup (like Tseytin did with $\fC$) and subsequently encode those relations into some fewer relations with basic encodings (like Makanin and Matiyasevich did); but in $\mathfrak{N}$, we must also heavily compress the relations into essentially a single relation $W_1 = W_2$ together with the first two auxiliary relations. In particular, new techniques for such compression had to be introduced. Note that Matiyasevich's original 1967 article \cite{Matiyasevich1967} does not contain any proof of undecidability of the word problem in \eqref{Eq:matiyasevich-semigroup}, only the statement of the results. The first detailed proof (other than that given in Matiyasevich's Ph.\ D.\ thesis, which is not readily available) using his original idea appeared only in 1995 in \cite{Matiyasevich1995}, although one can recover the idea from \cite{Boone1971}, and a somewhat different proof can be extracted from \cite{Collins1969}. The proof given in \cite{Matiyasevich1995} is, by contrast, self-contained (modulo the Markov--Post result of the undecidability of the word problem in semigroups) and does not, in particular, require the undecidability of the word problem in groups, unlike the path via the Tseytin semigroup $\fC$. We cannot here summarize the proof, and the reader interested in understanding \eqref{Eq:matiyasevich-semigroup} further is therefore directed to \cite{Matiyasevich1995}. This also contains material on and references about the termination problem for \textit{rewriting systems} (see \S\ref{Subsec:termination-problem}) with few rules. 

The word problem for monoids remains open both in the case of two and one relations. In recent years, there have been some developments on the subject of one-relation monoids. A more complete story of this problem and its history can be found in the author's recent (2021) survey of this problem \cite{NybergBroddaSpecial}. However, since then, there have been some new developments on the subject, which bear mentioning here. First, a one-relation monoid $\pres{Mon}{A}{u=v}$ will be said to be \textit{monadic} if its relation is of the form $u = a$, where $a \in A$ is a single letter (and $u$ is a non-empty word). It is unknown (!) whether the word problem is always decidable for monadic one-relation monoids. By contrast, for special monoids of the form $\pres{Mon}{A}{w=1}$ the word problem is known to be decidable, by the Adian--Makanin Theorem in \S\ref{Subsec:specialmonoids} combined with the decidability of the word problem in one-relator groups, as proved by Magnus \cite{Magnus1932}). Oganesian \cite{Oganesian1978, Oganesian1982} studied monadic one-relation monoids in detail. We refer the reader to \cite[\S6.3]{NybergBroddaSurvey} for a more complete account. 

A remarkable contribution to the theory of monadic one-relation monoids was given by Guba \cite{Guba1997} who in 1997 proved several interesting theorems; one of them linked \textit{membership problems} for one-relator groups to the word problem for monadic one-relation monoids. Motivated in part by this link, the author, Foniqi \& Gray \cite{IsoGrayCF} have demonstrated the general undecidability of the \textit{rational subset membership problem} in monadic one-relation monoids. This represents the first undecidable (non-artificial) problem for one-relation monoids, and may indicate that the record of three relations for an undecidable word problem may one day be reduced to just one relation. No matter the outcome of this problem, combinatorial semigroup theory still has much to offer in way of undecidability results.

\subsection{Termination for rewriting systems}\label{Subsec:termination-problem} We shall now briefly forego the chronological order of the consequences of Tseytin's result in order to stay on topic; we will return to the chronology in \S\ref{Subsec:group-consequence-borisov}. For now, instead, we will present a topic closely related to semigroups, for which the word problem is perhaps an even more natural problem. 

A \textit{rewriting system} may be thought of as a ``one-sided semigroup presentation''; a semisemigroup of sorts, in which one may only apply the relations of a presentation from one direction to the other. Such systems carry other names in the literature, where perhaps the most common is \textit{semi-Thue systems}, after A. Thue, who studied some simple one-rule systems in \cite{Thue1914}. We shall here retain, for simplicity, the notation \textit{rewriting system}. Formally, a rewriting system over an alphabet $A$ is a subset $\mathcal{S} \subseteq A^\ast \times A^\ast$, and is hence syntactically no different from a monoid presentation. However, we will normally always be favouring a semantic interpretation of these subsets; first, if $(u, v) \in \mathcal{S}$ then we will write this as $(u \to v) \in \mathcal{S}$, and call this a \textit{rule} of $\mathcal{S}$. A simple example comes from the one-sided version of the monoid defined in \eqref{Eq:monoid-ab-ba-pres}, being the rewriting system $\mathcal{S} = \{ (ab,ba) \} \subseteq \{ a, b \}^\ast \times \{ a, b \}^\ast$, which we write more readably as
\begin{equation}\label{Eq:ab-ba-rewriting-example}
\text{Alphabet: } \{a, b\}; \quad \text{Rule(s): } ab \to ba.
\end{equation}
Importantly, we will also induce several relations on $A^\ast$ from $\mathcal{S}$. We will think of the rules of $\mathcal{S}$ as being applicable to words over the alphabet of $\mathcal{S}$. If $w_1 = xuy$ and $w_2 = xvy$, where $x, y \in A^\ast$ and $u \to v$ is a rule of $\mathcal{S}$, then we say that $w_1$ \textit{rewrites} (in a single step) to $w_2$. Write this as $w_1 \xr{\mathcal{S}} w_2$. The transitive, reflexive closure of this operator is denoted $\xra{\mathcal{S}}$, and if $w_1 \xra{\mathcal{S}} w_2$ then we say that $w_2$ is a \textit{descendant} of $w_1$ (mod $\mathcal{S}$). If there is no infinite chain 
\[
w_1 \xr{\mathcal{S}} w_2 \xr{\mathcal{S}} w_3 \xr{\mathcal{S}} \cdots 
\]
(where the $w_i$ are not necessarily distinct) then we say that $\mathcal{S}$ is \textit{terminating} (or \textit{Noetherian}). Thus $\xra{S}$ is the one-sided analogue of $=_S$ defined in \S\ref{Sec:background}.

A word is said to be \textit{irreducible} if no rule of $\mathcal{S}$ can be applied to it. Given a terminating rewriting system, every element $w \in A^\ast$ can be rewritten to some irreducible descendant by successively applying rules; however, there may be many distinct irreducible descendants of a given word. Any system in which every word has a unique irreducible descendant is called \textit{uniquely terminating} or \textit{complete}. An equivalent condition (assuming termination) is that of \textit{confluence}: for every $u, w_1, w_2 \in A^\ast$, if $u \xra{\mathcal{S}} w_1$ and $u \xra{\mathcal{S}} w_2$, then there exists some $v \in A^\ast$ such that $w_1 \xra{\mathcal{S}} v$ and $w_2 \xra{\mathcal{S}} v$. A system is \textit{locally confluent} if the following somewhat weaker condition holds, being only concerned with direct descendants: for every $u, w_1, w_2 \in A^\ast$, if $u \xr{\mathcal{S}} w_1$ and $u \xr{\mathcal{S}} w_2$, then there exists some $v \in A^\ast$ such that $w_1 \xra{\mathcal{S}} v$ and $w_2 \xra{\mathcal{S}} v$. Newman's Lemma \cite{Newman1942} shows that for a terminating system, local confluence is equivalent to confluence. In particular, given a system with a single rule $(w \to 1)$, where $w$ has no non-trivial suffix which is also a prefix (i.e.\ $w$ is \textit{self-overlap free}), then this system is confluent. For example, the system with the single rule $(bc \to 1)$ is confluent. Furthermore, the system \eqref{Eq:ab-ba-rewriting-example} is complete.

Any rewriting system $\mathcal{S}$ gives rise to an associated monoid $M(\mathcal{S})$, the latter simply being the monoid with the same generators as $\mathcal{S}$, and a defining relation $u=v$ whenever $u \to v$ is a rule of $\mathcal{S}$. If $M$ is a monoid such that $M = M(\mathcal{S})$ for some rewriting system $\mathcal{S}$, then we say that $M$ \textit{admits} the rewriting system $\mathcal{S}$. For example, the free commutative monoid $\mathbb{N}^2$ admits the finite complete rewriting system \eqref{Eq:ab-ba-rewriting-example}. More generally, consider the one-relation monoid $M = \pres{Mon}{A}{u=v}$ with $|u| \geq |v|$ and where $u$ is self-overlap free. Then $u \to v$ is a complete rewriting system for $M$. However, if we drop the condition that $|u| \geq |v|$ then almost nothing is known about one-relation monoids admitting finite complete rewriting systems; indeed, as in \S\ref{Subsec:three-relations}, the word problem remains open for one-relation monoids of the form $\pres{Mon}{a,b}{a=v}$, where $v$ is arbitrary. 

It is clear that any monoid $M$ admitting a finite complete rewriting system $\mathcal{S}$ has decidable word problem: given two words $u, v$, simply rewrite them using $\mathcal{S}$ to their irreducible descendants, and compare those. We have $u = v$ in $M$ if and only if the two descendants are identical as words.\footnote{The time complexity of this solution can, however, be arbitrarily worse than that of the word problem for $M$, see \cite{Bauer1984}.} Thus, it is natural to ask whether every monoid with decidable word problem admits a finite complete rewriting system. It is at least conceivable to imagine how one might give a positive answer to such a question: namely, to construct such a rewriting system for any given monoid with decidable word problem. But how to prove a negative answer is not so clear: after all, how could one prove that a given monoid does not admit \textit{any} finite complete rewriting system? After all, monoids admit many presentations; how can one prove that \textit{none} of them are of a particular form? The missing tool to do so comes via homological algebra, and while this is entirely outside the scope of this article, we do mention the resolution of the question: there exists a monoid $M$ with decidable word problem, but which does not admit any finite complete rewriting system. The key comes from the fact that any monoid $M$ admitting a finite complete rewriting system must also satisfy certain homological finiteness properties; Squier \cite{Squier1987} proved that such $M$ must admit $\operatorname{FP}_3$, and Kobayashi \cite{Kobayashi1990} extended this to $\operatorname{FP}_\infty$. But there are many monoids -- even groups, as constructed by Stallings \cite{Stallings1963} -- with decidable word problem and which nevertheless do not even satisfy $\operatorname{FP}_3$. This yields the negative answer to the question. We note also that while the problem of whether all one-relation monoids admit finite complete rewriting systems remains open, Gray \& Steinberg \cite{Gray2022Steinberg} have recently shown that all one-relation monoids are $\operatorname{FP}_\infty$. 

How does rewriting systems relate to Tseytin's semigroup $\fC$? As before, the encoding ideas introduced by Tseytin have been used very fruitfully here. The key problem here is the \textit{termination problem}. For a rewriting system $\mathcal{S} \subseteq A^\ast \times A^\ast$, this is the problem of deciding, upon given a word $w \in A^\ast$, whether or not there exists an infinite sequence $w \xr{\mathcal{S}} w_1 \xr{\mathcal{S}} w_2 \xr{\mathcal{S}} \cdots$ starting at $w$. If there is no such derivation, then we say that $\mathcal{S}$ \textit{terminates} on $w$.  Unsurprisingly, the termination problem is undecidable in general for finite systems $\mathcal{S}$. This undecidability can be accomplished by encoding the behaviour of a universal Turing machine into the rules of $\mathcal{S}$, and the behaviour of any given Turing machine $T$ into a word $w = w(T)$, in such a way that $\mathcal{S}$ terminates on $w(T)$ if and only if $T$ halts (on the empty input, say). This was first done by Huet \& Lankford \cite{Huet1978}. Such encodings will, however, generally have very many rules. Thus, in pursuing the above question, it is of interest to understand how few rules an $\mathcal{S}$ with undecidable termination problem can have. Indeed, the following problem remains tantalizingly open:

\begin{question*}
Is the termination problem decidable for all one-rule rewriting systems? 
\end{question*}

This problem was probably first posed explicitly in \cite[p.\ 110]{Dauchet1989}. S\'{e}nizergues \cite[Theorem~19]{Senizergues1995} was able to construct a $10$-rule rewriting system $\mathcal{S}$ such that the termination problem for $\mathcal{S}$ is undecidable. This was improved by Matiyasevich \& S\'{e}nizergues \cite{Matiyasevich2005} to a $3$-rule rewriting system with undecidable termination problem, and the proof uses the $3$-relation semigroup \eqref{Eq:matiyasevich-semigroup} of Matiyasevich with undecidable word problem. However, the termination result is \textit{not} a direct consequence of the undecidability of the word problem in \eqref{Eq:matiyasevich-semigroup}. Indeed, the article \cite{Matiyasevich2005} introduces a wealth of new ideas in this context for proving undecidability of the termination problem. The resulting $3$-rule rewriting system with undecidable termination problem shares one drawback of \eqref{Eq:matiyasevich-semigroup}: one of its rules has both sides be extremely long words. It was therefore of interest to construct an explicit example of simpler example, even if it required more rules. 

A recent article by Halava, Matiyasevich \& Niskanen \cite{Halava2017} accomplishes precisely this. What is more, it uses, as explicitly acknowledged in their abstract, ideas by Tseytin. Their system $\mathcal{U}$ has $24$ rules, and each of the rules $(u \to v)$ has both $u$ and $v$ have length $\leq 5$. Using an encoding function $\varphi$ (entirely analogous to the kind described in \eqref{Eq:M1-dec-M2-monoid}), they prove that $\mathcal{U}$ is universal in the following sense: for any given rewriting system $\mathcal{S} \subseteq A^\ast \times A^\ast$, and for any $w \in A^\ast$, we have that $\mathcal{S}$ terminates on $w$ if and only if $\mathcal{U}$ terminates on the word $\varphi(w)$. Thus, since the termination problem is undecidable for \textit{some} rewriting system, it follows that the termination problem is undecidable for $\mathcal{U}$. 

We will not write the rules of the rewriting system $\mathcal{U}$, but we will say that there is a large shadow cast on them by the relations of Tseytin's semigroup $\fC$. Instead of one ``stylus'' letter $t$ (as described in \S\ref{Sec:Tseytin-proof}), the system $\mathcal{U}$ has four such styluses $e_1, e_2, e_3, e_4$. If one disregards the role played by them, it does not take long to see the relations of $\fC$ appear in the rules of $\mathcal{U}$. The proof of the aforementioned universality of $\mathcal{U}$ follows a similar pattern as the proof of the universality of $\fC$, too; however, several new ideas, unique to rewriting systems, are required. 

The $3$-rule system of Matiyasevich \& S\'{e}nizergues remains the smallest, with respect to the number of rules, rewriting system with undecidable termination problem. The termination problem for one-rule rewriting systems, like the word problem for one-relation monoids, remains a tantalizingly difficult problem. Indeed, the termination problem represents a strikingly fundamental question in theoretical computer science, as it is indirectly asking for an answer to the question: \textit{when can we decide whether a given computation will terminate or not?} The fact that we do not know whether we can solve this question when said computation uses only a single rule (in this particular model of computation) is frustrating; the fact that the aforementioned contributions to the problem all stem from work by Tseytin clearly demonstrates the continued importance of his contributions to this area, and by extension all theoretical computer science. 

\subsection{The word problem for groups}\label{Subsec:group-consequence-borisov} There are many beautiful consequences of Tseytin's result, and in the previous sections we gave a brief overview of its impact on semigroups and rewriting systems (i.e.\ semisemigroups). In this section, we will present perhaps the most remarkable endpoint of the series of consequences of Tseytin's result, especially as it relates to less specialized areas than semigroup theory. This is the 1969 construction by Borisov \cite{Borisov1969} of a 12-relator group with undecidable word problem. 

We shall not go in-depth with this construction, but we will give a very brief overview of the idea. The modern reader familiar with combinatorial group theory will recognise many familiar statements in Borisov's article, if in a different language; there is the notion of an HNN-extension (via a \textit{base} and \textit{filtering letters}, today called the \textit{base group} and \textit{stable letters}), and Britton's lemma, both in its embedding form as \cite[Lemma~1]{Borisov1969} and its form as a statement about pinched subwords as \cite[Lemma~2]{Borisov1969}. The modern reader not familiar with combinatorial group theory is warmly recommended to consult \cite{Baumslag1993,Lyndon1977,Magnus1966} for excellent introductions. 

Returning to Borisov's article, the key idea is to use combinatorial techniques in groups -- specifically HNN-extensions and Britton's lemma -- to ``simulate'' the very (un)pleasant combinatorial nature of semigroups. Let us first give an example of this nature in semigroups, and how such na\"ive constructions fail in groups. For example, consider the monoid with three generators $a, b, c$ and defining relations $ab^{i^2}c = ab^{j^2}c$ for all $i, j \in \N$. Which words are equal, in this monoid, to $ab^2c$? Since $2$ is not a perfect square, only $ab^2c$ is equal to $ab^2c$ -- no relation can be applied to this word! On the other hand, $ab^nc = ab^4c$ in this monoid if and only if $n$ is a perfect square. Similarly, we may simulate cubes, primes, or any other number-theoretic construction using such monoids. What if we were to try something analogous in groups? Immediately, one runs into problems: in any group, the relation $ab^ic = ab^jc$ is equivalent to $b^i = b^j$. Hence, the structure of the group will collapse, and all number-theoretic subtlety is lost.

One might believe that such combinatorial tricks, even when they do work for semigroups, are fairly meaningless exercises in algebraic manipulations. But this is not so. Indeed, a trick not entirely unlike one of the above form is used by Markov \cite{Markov1951} to prove the undecidability of the isomorphism, triviality, and other problems for semigroups. In fact, such tricks are also used to prove the undecidability of Markov properties for \textit{groups}, as in the Adian--Rabin theorem \cite{Adian1955,Rabin1958}. To accomplish this, one must instead simulate the combinatorial nature of the tricks via HNN-extensions rather than na\"ively into the defining relations. This requires significantly more intricate effort to materialise. We emphasize that this simulatory nature of HNN-extensions is already explicitly present in the proof by Britton of the undecidability of the word problem in groups \cite{Britton1958}; this article is also cited by Borisov. Indeed, one might for this, and other, reasons argue that Britton's 1958 article \cite{Britton1958} ought to be considered the first article in combinatorial group theory. But this is a digression from the task at hand.  

With these simulatory ideas in mind, then, the broad outline of the article is easy to follow. Borisov starts with Matiyasevich's $3$-relation semigroup \eqref{Eq:matiyasevich-semigroup}, and from it produces a $14$-relator group $\Gamma$ which can simulate the word problem of this semigroup. Hence, $\Gamma$ has undecidable word problem. He then proceeds to embed $\Gamma$ into a $12$-relator group $\Gamma'$; hence the word problem is undecidable also in $12$-relator groups. But his method of proof is more general than this: it is a general construction, applicable to any semigroup with $n$ generators and $r$ defining relations, and which produces a group with $2n+r+5$ relators and undecidable word problem. As discussed in \S\ref{Sec:background}, one may always, in this context, assume $n=2$, so this shows that if the word problem is undecidable in some $r$-relation semigroup, then there is a $(9+r)$-relator group with undecidable word problem. Taking Matiyasevich's example, having $r=3$, then yields the desired $12$-relator group. Furthermore, a curious corollary of Borisov's result is thus: \textit{if the word problem is decidable in every ten-relator group, then the word problem is decidable in every one-relation monoid}. Thus a semigroup-theoretic problem is reduced to a group-theoretic problem. A word of caution is required. The word problem is known to be decidable in all one-relator groups; this is a classical result by Magnus \cite{Magnus1932}. But already the case of the word problem for two-relator groups is wide open (see \cite[Problem~9.29]{Kourovka2023}), so there seems exceedingly little hope that the word problem will be resolved positively for ten-relator groups anytime soon.

Borisov's arguments, particularly the step from the $14$-relator group to the $12$-relator group, have been simplified by Collins \cite{Collins1972}. Collins \cite{Collins1969} had also given a proof in 1969 of the existence of a $14$-relator group with undecidable word problem at the same time as Borisov. The argument in \cite{Collins1969} is also used to produce an $11$-relator group with undecidable conjugacy problem. 

A final remark is in order. The first page \cite[p. 521]{Borisov1969}\footnote{In the English translation in \textit{Math. Notes.} \textbf{6}, 1969, pp. 768--775, this is p.~768.} of Borisov's article contains a curious note. First, it gives a timeline of producing groups with few defining relations and undecidable word problem, and it is mentioned that using Matiyasevich's three-relation semigroup, one can push through the construction used by Boone \cite{Boone1958} to find a 28-relator group with undecidable word problem. Borisov then, remarkably, notes that Tseytin has produced a $14$-relator group with undecidable word problem, but that no proof had appeared at the time. Even today, no proof has appeared of this result by Tseytin, and I know of no other reference to it in the literature. It would be curious to see whether Tseytin's idea, assuming that it was sound, was similar to that of Borisov, or indeed to that of Collins. Given the generally rather convergent nature of the techniques of this subject, it would not be unreasonable to believe so. 

\subsection{Universal (semi)groups}\label{Subsec:modern-proof-Matiyasevich}

As with any good idea, Tseytin's construction has filtered through a series of developments and has been reduced to a core. We have seen many reflections of this core throughout this present article. The most modern form of this core appears in an article by Boone, Collins \& Matiyasevich \cite{Boone1971}. In this article, the authors prove the following remarkable theorem:

\begin{theorem*}[Boone, Collins \& Matiyasevich, 1971]
There exists a fixed six-relation semigroup containing an isomorphic copy of every finitely presented semigroup. 
\end{theorem*}

The proof goes as follows. First, given a finitely presented semigroup $S$, one constructs a semigroup $\mathfrak{D}$ which has as generators the generators of $S$ together with new ``instruction'' generators, and as defining relations a (suitably encoded) copy of the relations of $S$, together with a system of ``instruction'' relations. These instruction relations are all written over the instruction generators, and in particular do not depend on the relations of $S$. Indeed, the instruction relations are very similar to the relations of $\fC$ (i.e.\ Tseytin's semigroup from \S\ref{Sec:Tseytin-proof}), but augmented somewhat to account for two styluses rather than one. It is then shown, using methods entirely analogous to those used by Tseytin, that two words are equal in $S$ if and only if they are equal in $\mathfrak{D}$. This shows that $S$ embeds isomorphically into $\mathfrak{D}$ (via the identity map). The number of defining relations of $\mathfrak{D}$ is $n+4$, where $n$ is the number of generators of $S$. By the result of Hall \cite{Hall1949} mentioned in \S\ref{Sec:background}, every semigroup can be embedded in a two-generator semigroup, and hence we find that every finitely presented semigroup can be embedded into some six-relation semigroup. To get the above theorem, we must now reverse the quantifiers. To do this, we finally appeal to a result by Murskii \cite{Murskii1967}, who proved that there exists a fixed \textit{universal} finitely presented semigroup $\mathfrak{M}$, in the sense that every finitely presented (indeed every recursively presented) semigroup can be embedded into $\mathfrak{M}$. Since $\mathfrak{M}$ can now be embedded in a six-relation semigroup, this yields the above theorem. We have omitted all technical details of the proof; excluding the result by Murskii, these are no more technical than Tseytin's. With this in mind, we refer the reader to \cite{Boone1971}.

Finally, the above sort of universality questions have also been studied for groups. The analogue of Murskii's result for groups is called \textit{Higman's embedding theorem} \cite{Higman1961}. Using this and similar types of embedding constructions as those of Boone, Collins \& Matiyasevich, Valiev \cite{Valiev1973} gave in 1973 the first explicit example of a universal finitely presented group, being a $42$-relator group containing every finitely presented group. This was improved by Boone \& Collins \cite{Boone1974} the next year to $26$ relators. The ideas contained therein are entirely analogous to the above semigroup-theoretic result. Finally, Valiev \cite{Valiev1977} improved this by constructing a 21-relator universal group, and this has not been improved since. Nevertheless, many natural questions of this type continue to be the subject of modern research. For example, in 2014, Chiodo \cite{Chiodo2014} proved that there is a universal \textit{torsion-free} finitely presented group, in the sense that there is a finitely presented group $H$ which is torsion-free, and with the property that $H$ contains an isomorphic copy of every finitely presented torsion-free group. Doubtlessly, many other results remain provable but yet unproved; and equally doubtlessly, Tseytin's ideas will continue to resonate through their eventual proofs.

\bibliography{TseytinNybergBrodda.bib}
\bibliographystyle{plain}
 


\section*{Supplementary material}

G. S. Tseytin

# AN ASSOCIATIVE CALCULUS WITH UNSOLVABLE EQUIVALENCE PROBLEM

## §1

In this work we will consider associative calculi, as defined in [2, VI. §1.1–4]. The equivalence problem for associative calculi will be understood in the sense of [2, VI. §1.9] (as a normal mass problem [2, V.]).

There exist associative calculi for which the problem of equivalence is unsolvable, which was first proved by A. A. Markov [1] and E. L. Post [6]. This result (as well as the existence of associative calculi in which the problem of equivalence to a fixed word is unsolvable [2, VI. §2.4, §2.3.2]) serves as the source for a number of theorems on the existence of unsolvable mass problems for other mathematical objects. For example, on the basis of this, in [2] is proved the existence of reversible Post calculi with unsolvable derivability problem [2, VI. §4.8.2], the existence of a system of integer matrices for which the representability problem is unsolvable [2, VI. §10]; in [3] is proved the existence of a group with unsolvable conjugacy problem and a two-dimensional polyhedron with unsolvable path homotopy problem [3, §6]. In constructing the corresponding mathematical objects (the Post calculus, the system of integer matrices, etc.) the original associative calculus is used. Therefore, to construct simpler objects with such properties, it is of interest to find a specific, sufficiently simple, associative calculus with an unsolvable equivalence problem (respectively a calculus with an unsolvable problem of equivalence to a fixed word). Additionally, such a calculus could be the source of new results regarding unsolvable mass problems. A. A. Markov in [2, VI. §6] constructed an associative calculus with an unsolvable equivalence problem with a 13-letter alphabet and with 33 defining relations. In the present work, we will give even simpler examples of such associative calculi. Namely, we will prove that the equivalence problem is unsolvable in the calculus $\mathfrak{C}_1$ on the five-letter alphabet $\{a, b, c, d, e\}$ and with the following defining relations:

$$(1) \quad \begin{cases} ac \longleftrightarrow ca \\ ad \longleftrightarrow da \\ bc \longleftrightarrow cb \\ bd \longleftrightarrow db \\ eca \longleftrightarrow ce \\ edb \longleftrightarrow de \\ cca \longleftrightarrow ccae. \end{cases}$$

Furthermore, we will prove that in the calculus $\mathfrak{C}_1$ the problem of equivalence with some fixed word is unsolvable, but this word is extremely long. We therefore construct a

calculus $\mathfrak{C}_2$ also on the alphabet $\{a, b, c, d, e\}$ and with the defining relations:

$$(2) \quad \begin{cases} ac \longleftrightarrow ca \\ ad \longleftrightarrow da \\ bc \longleftrightarrow cb \\ bd \longleftrightarrow db \\ eca \longleftrightarrow ce \\ edb \longleftrightarrow de \\ cdca \longleftrightarrow cdcae \\ caaa \longleftrightarrow aaa \\ daaa \longleftrightarrow aaa. \end{cases}$$

In this calculus, we will prove that the problem of equivalence with the word $aaa$ is unsolvable.

All our unsolvability results will be based on the proof of P. S. Novikov [4] of the unsolvability of the identity problem in groups. In doing this, we will have to reason about some normal mass problems and prove the reducibility of some problems to others. In order to give an exact formulation of these normal mass problems, it is necessary to fix the way of writing the input data over some alphabet, and then pose the problem of finding some normal algorithm over this alphabet. For simplicity, we will not do this, and instead implicitly imply that we are speaking of normal mass problems obtained by a proper refinement of this formulation.[1] The principle of normalization [2, II. §5] allows us to state that such a refinement of the mass problem in question is uniquely defined up to equivalence. In the same way, proofs of the reducibility of some normal mass problem to another, which formally requires reasoning using normal algorithms, will be replaced, again for simplicity, by a description of the processes by which the solution of a single problem corresponding to one mass problem can be corresponded to one of another mass problem. An absolutely precise presentation of all the proofs of this kind given here poses no fundamental difficulty except its great cumbersomeness.

This article uses terms and notation from [2], in particular the notation for letter substitution and the projection operations [2, II. §4.14]. The article is divided into sections, within which lemmas are separately numbered (by numbers with dots), as are formulae (by numbers in parentheses). A reference to a lemma or a formula is made by indicating the section number and the number of the lemma or formula; the entire reference is enclosed in square brackets. For references within the same section, the section number is omitted.[2]

An abbreviated presentation of the results in this article appear (in a slightly modified form) already in [5].

---

[1] Some difficulty in refining the formulations of mass problems is posed by those problems that in their original formulation use the concept of an arbitrary alphabet. This difficulty can be overcome, for example, by encoding words [2, I. §6.2] into some fixed alphabet and using [2, Theorem 6.2.8].

[2] Any reference made to articles in the bibliography at the end of the article is given in square brackets by a number without a dot.

## §2

We introduce some new operations for encoding words (in the sequel, the term "encoding" will **not** be used in the same sense as in [2]). Consider an arbitrary alphabet $\Pi$ and a two-letter alphabet $A_0 = \{a, b\}$. Let us associate to each letter from $\Pi$ some non-negative integer such that different letters are associated with different numbers; the number associated to the letter $\xi$ will be denoted by $n(\xi)$. If $\xi \neq \eta$, then $n(\xi) \neq n(\eta)$. We let the *left encoding* of a letter $\xi$ from $\Pi$ be the word $ab^{n(\xi)}$ over $A_0$, and the *right encoding* of the letter $\xi$ will be the word $b^{n(\xi)}a$. We let the *left encoding* of a word $P$ in $\Pi$ be the word $[P^\zeta$ obtained from $P$ by replacing each of its letters by its left encoding; the *right encoding* of $P$ will be the word $[P^\supset$ obtained from $P$ by replacing each of the letters of $P$ by its right encoding.

(1) $$[P^\zeta = \mathbf{S}^{\cdots\xi\cdots}_{\cdots ab^{n(\xi)}\cdots}P|.$$

(2) $$[P^\supset = \mathbf{S}^{\cdots\xi\cdots}_{\cdots b^{n(\xi)}a\cdots}P|.$$

Both of $[P^\zeta$ and $[P^\supset$ are words over $A_0$. Clearly, for any words $P$ and $Q$ in $\Pi$, we have

(3) $$[PQ^\zeta = [P^\zeta[Q^\zeta,$$

(4) $$[PQ^\supset = [P^\supset[Q^\supset.$$

**Lemma 1.** *For any word $P$ in $\Pi$,*

(5) $$a[P^\zeta = [P^\zeta a.$$

This statement is proved without difficulty using right induction on $\Pi$, see [2, I.§3.8]. The *encoding* of the word $P$ over $\Pi$ will be the word $[P^\bigcirc$ defined as $a[P^\supset$.

(6) $$[P^\bigcirc = a[P^\supset = [P^\zeta a. \qquad\qquad [(5)$$

Let $P, Q$, and $R$ be words over $\Pi$. Then

$$[PQR^\bigcirc = a[PQR^\supset \qquad\qquad [(6)$$
$$= a[P^\supset[Q^\supset[R^\supset \qquad\qquad [(4)$$
$$= [P^\zeta a[Q^\supset[R^\supset \qquad\qquad [(5)$$
$$= [P^\zeta[Q^\bigcirc[R^\supset. \qquad\qquad [(6)$$

Hence, for all words $P, Q$, and $R$ over $\Pi$, we have

(7) $$[PQR^\bigcirc = [P^\zeta[Q^\bigcirc[R^\supset,$$

(8) $$[PQ^\bigcirc = [P^\zeta[Q^\bigcirc, \qquad\qquad [(7)$$

(9) $$[QR^\bigcirc = [Q^\bigcirc[R^\supset \qquad\qquad [(7)$$

**Lemma 2.** *If $P$ and $Q$ are words over $\Pi$, and $[P^\bigcirc = [Q^\bigcirc$, then $P = Q$.*

*Proof.* Let $R$ be the longest common prefix of $P$ and $Q$ (see [2, I.§3.11]). Then there exists some left coprime words $S$ and $T$ such that $P = RS$ and $Q = RT$. By [(9)], we have $[P^\bigcirc = [R^\bigcirc[S^\supset$, and $[Q^\bigcirc = [R^\bigcirc[T^\supset$. As $[P^\bigcirc = [Q^\bigcirc$, it follows that $[S^\supset = [T^\supset$. If $S \neq \Lambda$, then $[S^\supset \neq \Lambda$, and $T \neq \Lambda$. Suppose the leftmost letter of $S$ is $\xi$, and that the leftmost letter of $T$ is $\eta$. Then $[S^\supset$ begins with the word $b^{n(\xi)}a$, and $[T^\supset$ begins with the word $b^{n(\eta)}$. But $[S^\supset = [T^\supset$, so either $b^{n(\xi)}a$ is a prefix of $b^{n(\eta)}a$, or else $b^{n(\eta)}a$ is a prefix of $b^{n(\xi)}a$. But this is possible only if $n(\xi) = n(\eta)$, and then $\xi = \eta$. Hence the first letters of $S$ and $T$ coincide, which contradicts the left coprimality of $S$ and $T$. Hence, $S$ cannot be non-empty, i.e. $S = \Lambda$. But then $[T^\supset = [S^\supset = \Lambda$, so $T = \Lambda$. Consequently, $P = R$ and $Q = R$, i.e. $P = Q$, which is what was to be shown. $\qquad\square$

**Lemma 3.** *If $P$ is a word over $\Pi$, $M$ and $N$ are words over $A_0$, and $\lfloor P^{\circ} = MaN$, then there exists some words $Q$ and $R$ over $\Pi$ such that $M = \lfloor Q^{\zeta}$ and $P = QR$.*

*Proof.* We will prove the result by induction on the length of the word $M$. When $M = \Lambda$, then $Q$ can be taken to be $\Lambda$, and $R$ can be taken to be $P$. Suppose that we have proved the claim whenever the length of $M$ is less than $m$ (with $m > 0$). Suppose then that the length of $M$ is $m$. Then $M$ is non-empty. Hence, $MaN = \lfloor P^{\circ} = a \lfloor P^{\gimel}$ by $\lfloor (6)\rfloor$, so $M$ begins with the letter $a$, i.e. it can be written as $aM_1$, where $M$ is some word over $A_0$. We have $\lfloor P^{\circ} \neq a$, as, by Lemma 2, $P \neq \Lambda$. Thus $P$ can be written as $\xi P_1$, where $\xi \in \Pi$ and $P_1$ is a word over $\Pi$. Then, by $\lfloor (8)\rfloor$, we have $\lfloor P^{\circ} = \lfloor \xi^{\zeta} \lfloor P_1{}^{\circ} = ab^{n(\xi)} \lfloor P_1{}^{\circ}$ and $ab^{n(\xi)} \lfloor P_1{}^{\circ} = aM_1 aN$. The word $ab^{n(\xi)}$ begins with the word $aM_1a$, so $aM_1a$ begins with the word $ab^{n(\xi)}$, i.e. there exists a non-empty word $M_2$ over $A_0$ such that $aM_1a = ab^{n(\xi)} M_2$. The word $M_2$ ends with the letter $a$, so we can write $M_2 = M_3 a$, where $M_3$ is a word over $A_0$. We have $aM_1 = ab^{n(\xi)} M_3$, so the length of $M_3$ is less than $m$. Furthermore, $ab^{n(\xi)} \lfloor P_1{}^{\circ} = aM_1 aN = ab^{n(\xi)} M_3 aN$, so $\lfloor P_1{}^{\circ} = M_3 aN$. By the inductive hypothesis, there exist some words $Q_1$ and $R$ in $\Pi$ such that $M_3 = \lfloor Q_1{}^{\zeta}$ and $P_1 = Q_1 R$. We have $Q = \xi Q_1$. Then $\lfloor Q^{\zeta} = \lfloor \xi^{\zeta} \lfloor Q_1{}^{\zeta} = ab^{n(\xi)} M_3 = M$ and $QR = \xi Q_1 R = \xi P_1 = P$. This completes the proof. □

**Lemma 4.** *If $P$ and $Q$ are words over $\Pi$, $M$ and $N$ are words over $A_0$, and $\lfloor P^{\circ} = M \lfloor Q^{\circ} N$, then there exist some words $R$ and $S$ over $\Pi$ such that $\lfloor R^{\zeta} = M$, $\lfloor S^{\gimel} = N$, and $P = RQS$.*

*Proof.* We have $\lfloor P^{\circ} = Ma \lfloor Q^{\gimel} N$, and so by Lemma 3 there exist some words $R$ and $T$ over $\Pi$ such that $P = RT$ and $\lfloor R^{\zeta} = M$. Then $\lfloor P^{\circ} = \lfloor R^{\zeta} \lfloor T^{\circ}$ by $\lfloor (8)\rfloor$, so $M \lfloor T^{\circ} = \lfloor P^{\circ} = M \lfloor Q^{\circ} N$ and $\lfloor T^{\circ} = \lfloor Q^{\circ} N = \lfloor Q^{\zeta} aN$. Again by Lemma 3 there exist words $Q_1$ and $S$ over $\Pi$ such that $T = Q_1 S$ and $\lfloor Q_1{}^{\zeta} = \lfloor Q^{\zeta}$. Then $\lfloor Q_1{}^{\circ} = \lfloor Q^{\circ}$, so by Lemma 2 we have $Q_1 = Q$. Since $T = QS$, we have $\lfloor T^{\circ} = \lfloor Q^{\circ} \lfloor S^{\circ}$ by $\lfloor (9)\rfloor$, so $\lfloor Q^{\circ} \lfloor S^{\gimel} = \lfloor Q^{\circ} N$ and $\lfloor S^{\gimel} = N$. Finally, as $P = RT$ and $T = QS$, we have $P = RQS$, completing the proof. □

Suppose $\mathfrak{P}$ is some associative calculus on the alphabet $\Pi$. We will define an associative calculus $\mathfrak{P}^{\circ}$ on the alphabet $A_0$, which will have defining system obtained from the defining system of $\mathfrak{P}$ by replacing all left and right sides of each defining relation by their encodings. That is, $M \longleftrightarrow N$ will be a defining relation of the system $\mathfrak{P}^{\circ}$ if and only if there exist words $P$ and $Q$ over $\Pi$ such that $\lfloor P^{\circ} = M$, $\lfloor Q^{\circ} = N$, and $P \longleftrightarrow Q$ is a defining relation of the system $\mathfrak{P}$. We will establish a link between the calculi $\mathfrak{P}$ and $\mathfrak{P}^{\circ}$.

**Lemma 5.** *If $P$ and $Q$ are words over $\Pi$, and $\mathfrak{P} : P \perp Q$, then $\mathfrak{P}^{\circ} : \lfloor P^{\circ} \perp \lfloor Q^{\circ}$.*

*Proof.* As $\mathfrak{P} : P \perp Q$, there exists some words $T, U, V$, and $W$ such that $P = VTW$, $Q = VUW$, and one of the relations $T \longleftrightarrow U$ and $U \longleftrightarrow T$ is a defining relation of the system $\mathfrak{P}$. But then, by $\lfloor (7)\rfloor$, we have that $\lfloor P^{\circ} = \lfloor V^{\zeta} \lfloor T^{\circ} \lfloor W^{\gimel}$ and $\lfloor W^{\circ} = \lfloor V^{\zeta} \lfloor U^{\circ} \lfloor W^{\gimel}$, and one of $\lfloor T^{\circ} \longleftrightarrow \lfloor U^{\circ}$ and $\lfloor U^{\circ} \longleftrightarrow \lfloor T^{\circ}$ will be a defining relation of the system $\mathfrak{P}^{\circ}$. Consequently, $\mathfrak{P}^{\circ} : \lfloor P^{\circ} \perp \lfloor Q^{\circ}$, which is what was to be shown. □

**Lemma 6.** *If $P$ and $Q$ are words over $\Pi$, and $\mathfrak{P} : P \perp\!\!\!\perp Q$, then $\mathfrak{P}^{\circ} : \lfloor P^{\circ} \perp\!\!\!\perp \lfloor Q^{\circ}$.*

This is easily proved by using Lemma 5 and induction on the number of entries in a $\mathfrak{P}$-sequence [2, VI. §1.6] connecting $P$ and $Q$.

**Lemma 7.** *If $P$ is a word over $\Pi$, $M$ is a word over $A_0$, and $\mathfrak{P}^\circ : \lceil P^\circ \perp M$, then there exists some word $Q$ over $\Pi$ such that $M = \lceil Q^\circ$ and $\mathfrak{P} : P \perp Q$.*

*Proof.* There exist some words $T, U, V, W$ over $A_0$ such that $\lceil P^\circ = VTW$, $M = VUW$, and one of the relations $T \longleftrightarrow U$ and $U \longleftrightarrow T$ is a defining relation of the system $\mathfrak{P}^\circ$. But then there must exist some words $K$ and $L$ over $\Pi$ such that $\lceil K^\circ = T$, $\lceil L^\circ = U$, and among the defining relations of $\mathfrak{P}$ there is one of $K \longleftrightarrow L$ and $L \longleftrightarrow K$. Hence, $\lceil P^\circ = V \lceil K^\circ W$, and by Lemma 4 there exist some words $X$ and $Y$ over $\Pi$ such that $\lceil X^\subset = V$, $\lceil Y^\supset = W$, and $P = XKY$. Setting $Q = XLY$, we have $\mathfrak{P} : P \perp Q$. Furthermore, $M = VUW = \lceil X^\subset \lceil L^\circ \lceil Y^\supset = \lceil Q^\circ$, i.e. $M = \lceil Q^\circ$, which is what was to be shown.                                    $\square$

**Lemma 8.** *If $P$ is a word over $\Pi$, and $M$ is a word over $A_0$ such that $\mathfrak{P}^\circ : \lceil P^\circ \perp\!\!\!\perp M$, then there exists some word $Q$ over $\Pi$ such that $M = \lceil Q^\circ$ and $\mathfrak{P} : P \perp\!\!\!\perp Q$.*

This is easily proved using Lemma 7 by induction on the number of entries in a $\mathfrak{P}^\circ$-sequence connecting $\lceil P^\circ$ and $M$.

**Lemma 9.** *If $P$ and $Q$ are words over $\Pi$ and $\mathfrak{P}^\circ : \lceil P^\circ \perp\!\!\!\perp \lceil Q^\circ$, then $\mathfrak{P} : P \perp\!\!\!\perp Q$.*

This follows immediately from Lemma 8 and Lemma 2.

## §3

In this section we will introduce some concepts and objects needed in the sequel.

Let us consider the alphabets $A_0 = \{a, b\}$, $\overline{A}_0 = \{c, d\}$, and $A_1 = A_0 \cup \overline{A}_0 \cup \{e\} = \{a, b, c, d, e\}$. We will say that the letter $c$ is the *twin* of the letter $a$, and that the letter $d$ is the twin of the letter $b$. The twin of the letter $\xi$ (for $\xi \in A_0$) will be denoted $\bar{\xi}$. We introduce an operation transforming any word $P$ over $A_0$ into a word $\lceil P^-$, obtained from $P$ by replacing each of its letters by its twin:

$$(1) \qquad\qquad \lceil P^- = \mathbf{S}_{c,d}^{a,b} P \rceil.$$

We will introduce the reverse operation, applicable to words over $\overline{A}_0$:

$$(2) \qquad\qquad \lceil P^\wedge = \mathbf{S}_{a,b}^{c,d} P \rceil.$$

Clearly, $\lceil P^-$ is a word over $\overline{A}_0$ (if $P$ is over $A_0$) and $\lceil Q^\wedge$ is a word over $A_0$ (if $Q$ is over $\overline{A}_0$). If $P$ is a word over $A_0$, then, clearly,

$$(3) \qquad\qquad \lceil\lceil P^{-\wedge} = P. \qquad\qquad [(1), (2)]$$

Fix some non-negative integer $i$, and construct the associative calculus $\mathfrak{C}$ over the alphabet $A_1$ and with the following defining relations (written in an abbreviated way):

$$(4) \qquad \begin{array}{ll} (4_1) & \left\{ \begin{array}{ll} \xi\eta \longleftrightarrow \eta\xi & (\xi \in A_0, \eta \in \overline{A}_0) \\ e\bar{\xi}\xi \longleftrightarrow \bar{\xi}e & (\xi \in A_0) \\ cd^i ca \longleftrightarrow cd^i cae & \end{array} \right. \\ (4_2) & \\ (4_3) & \end{array}$$

For $i = 0$, this calculus is just the calculus $\mathfrak{C}_1$ defined in §1(1).

We will prove the unsolvability of the equivalence problem in $\mathfrak{C}$. For this we will establish a connection between the equivalence problem in $\mathfrak{C}$ and the equivalence problem in some other associative calculi.

We will call a *special associative calculus* any associative calculus whose defining relations all have the empty word as their right-hand side. Let $\mathfrak{A}$ be an arbitrary special

associative calculus, and let $A$ be the alphabet of $\mathfrak{A}$. The system of defining relations of $\mathfrak{A}$ will have the form

$$(5) \qquad \begin{cases} K_1 \longleftrightarrow \\ \cdots \\ K_m \longleftrightarrow, \end{cases}$$

where $K_1, \ldots, K_m$ are words over $A$ (and $m \geq 0$). Let $\alpha$ be a letter which does not appear in $A$. Let $B = A \cup \{\alpha\}$ and consider the associative calculus $\mathfrak{B}$ over $B$ with defining system

$$(6) \qquad \begin{cases} K_0\alpha \longleftrightarrow \\ K_1\alpha \longleftrightarrow \\ \cdots \\ K_m\alpha \longleftrightarrow, \end{cases}$$

where $K_0 = \Lambda$. The first of the relations [(6)] is $\alpha \longleftrightarrow$, and so

$$(7) \qquad \mathfrak{B} : \alpha \perp\!\!\!\perp \Lambda.$$

From this, it follows that for $1 \leq t \leq m$ we have $\mathfrak{B} : K_t \perp\!\!\!\perp K_t\alpha\perp\Lambda$, i.e.

$$(8) \qquad \mathfrak{B} : K_t \perp\!\!\!\perp \Lambda.$$

It is now not difficult to prove the following results.

**Lemma 1.** *If $P, Q$ are words over $A$ and $\mathfrak{A} : P\perp Q$, then $\mathfrak{B} : P \perp\!\!\!\perp Q$.*

The proof follows by using [(8)].

**Lemma 2.** *If $P, Q$ are words over $A$ and $\mathfrak{A} : P\perp Q$, then $\mathfrak{B} : P \perp\!\!\!\perp Q$.*

The proof follows by using Lemma 1.

**Lemma 3.** *If $P, Q$ are words over $A$ and $\mathfrak{B} : P\perp Q$, then $\mathfrak{A} : [P^A\perp[Q^A$ or $[P^A = [Q^A$.*

**Lemma 4.** *If $P, Q$ are words over $A$ and $\mathfrak{B} : P\perp Q$, then $\mathfrak{A} : [P^A \perp\!\!\!\perp [Q^A$.*

This is proved using Lemma 3.

We will now define an encoding from the alphabet $B$ into the alphabet $A_0$. Thus, the alphabet $B$ will play the role of the alphabet $\Pi$ from §2. The number $n(\xi)$ (for $\xi \in B$) will be chosen such that $n(\alpha) = i$. Then

$$(9) \qquad [\alpha^\varsigma = ab^i; \quad [\alpha^\circ = ab^ia. \qquad\qquad [\S2(1), \S2(6)]$$

Let us construct from $\mathfrak{B}$ a calculus $\mathfrak{B}^\circ$ in the same manner which $\mathfrak{P}^\circ$ was constructed from $\mathfrak{P}$. Then $\mathfrak{B}^\circ$ is a calculus over $A_0$ with defining system of relations

$$(10) \qquad \begin{cases} [K_0\alpha^\circ \longleftrightarrow a \\ [K_1\alpha^\circ \longleftrightarrow a \\ \cdots \\ [K_m\alpha^\circ \longleftrightarrow a. \end{cases}$$

By using Lemmas §2.6, §2.8, and §2.9, we can prove the following propositions.

**Lemma 5.** *If $P$ and $Q$ are words over $B$ and $\mathfrak{B} : P \perp\!\!\!\perp Q$, then $\mathfrak{B}^\circ : [P^\circ \perp\!\!\!\perp [Q^\circ$.*

**Lemma 6.** *If $P$ is a word over $B$, and $M$ is a word over $A_0$ such that $\mathfrak{B}^\circ : [P^\circ \perp\!\!\!\perp M$, then there exists some word $Q$ over $B$, such that $M = [Q^\circ$ and $\mathfrak{B} : P \perp\!\!\!\perp Q$.*

**Lemma 7.** *If $P$ and $Q$ are words over $B$ such that $\mathfrak{B}^{\circ} : [P^{\circ} \perp\!\!\!\perp [Q^{\circ}, \text{ then } \mathfrak{B} : P \perp\!\!\!\perp Q.$*

Consider the word

(11) $$\mathbf{S} = [[\alpha K_0 \alpha K_1 \alpha \cdots K_m \alpha^{\circ -}.$$

The word $\mathbf{S}$ is a word over $\overline{A}_0$. It is a kind of record of the defining relations of the system $\mathfrak{A}$. With the help of this record we will be able to prove the following link between the calculi $\mathfrak{C}$ and $\mathfrak{B}^{\circ}$: *if $P$ and $Q$ are words over $A_0$, then $\mathfrak{C} : \mathbf{S}P \perp\!\!\!\perp \mathbf{S}Q$ if and only if $\mathfrak{B}^{\circ} : P \perp\!\!\!\perp Q.$* The proof of this will be the objective of §4 and §5.

## §4

**Lemma 1.** *If $P$ is a word over $A_0$, $\xi \in \overline{A}_0$, and $M$ and $N$ are words over $A_1$, then*

$$\mathfrak{C} : MP\xi N \perp\!\!\!\perp M\xi PN.$$

*Proof.* This is easily proved by right induction on the word $P$ over $A_0$. Indeed, if $P$ is empty it is obvious, and if it is true for some $P$ and $\eta \in A_0$, then $[\![(4_1)]\!]$ will give

$$\mathfrak{C} : MP\eta\xi N \perp MP\xi\eta N \perp\!\!\!\perp M\xi P\eta N,$$

which proves the statement also for $P\eta$.    $\square$

**Lemma 2.** *If $P$ is a word over $A_0$, $Q$ a word over $\overline{A}_0$, and $M, N$ are words over $A_1$, then*

$$\mathfrak{C} : MPQN \perp\!\!\!\perp MQPN.$$

The proof follows by using Lemma 1 and right induction on $Q$ over $\overline{A}_0$.

**Lemma 3.** *If $P$ is a word over $A_0$, and $M$ and $N$ are words over $A_1$, then*

$$\mathfrak{C} : M[P^- eN \perp\!\!\!\perp Me[P^- PN.$$

*Proof.* To prove this, we use right induction on $P$ over $A_0$. If $P$ is empty, the statement is obviously true. Suppose it is proved for some $P$ over $A_0$ for any $M$ and $N$ over $A_1$, and let $\xi \in A_0$. Then

$$M[P\xi^- eN = M[P^- \bar{\xi} eN,$$
$$\mathfrak{C} : M[P^- \bar{\xi} eN \perp M[P^- e\bar{\xi}\xi N \qquad\qquad [\S 3(4_2)]$$
$$\perp\!\!\!\perp Me[P^- P\bar{\xi}\xi N \quad \text{(by the inductive hypothesis)}$$
$$\perp\!\!\!\perp Me[P^- \bar{\xi} P\xi N, \qquad\qquad\qquad [\text{Lemma 1}]$$
$$Me[P^- \bar{\xi} P\xi N = Me[P\xi^- P\xi N,$$

and finally $\mathfrak{C} : M[P\xi^- eN \perp\!\!\!\perp Me[P\xi^- P\xi N.$    $\square$

**Lemma 4.** *If $P$ and $Q$ are words over $A_0$ and $0 \le t \le m$ (where $t$ is an integer), then*

$$\mathfrak{C} : \mathbf{S}PaQ \perp\!\!\!\perp \mathbf{S}P[K_t \alpha^{\circ} Q.$$

*Proof.* The word $\alpha K_0 \alpha K_1 \alpha \cdots K_m \alpha$ contains an occurrence of the word $\alpha K_t \alpha$, as by $[\S 2(7)]$ the word $[\alpha K_t \alpha^{\circ}$ occurs in $[\alpha K_0 \alpha K_1 \alpha \cdots K_m \alpha^{\circ}$. Consequently by $[\S 3(11)]$ the

word $[[\alpha K_t \alpha^{\circ-}$ occurs in $\mathbf{S}$, i.e. $\mathbf{S}$ can be written as $V[[\alpha K_t \alpha^{\circ-} W$, where $V$ and $W$ are some words over $\overline{A}_0$. We have:

$$\mathbf{S}PaQ = V[[\alpha K_t \alpha^{\circ-} W PaQ,$$

$$\mathfrak{C} : V[[\alpha K_t \alpha^{\circ-} W PaQ \perp\!\!\!\perp V[[\alpha K_t \alpha^{\circ-} PaWQ \qquad \text{[Lemma 2]}$$

$$\perp\!\!\!\perp VP[[\alpha K_t \alpha^{\circ-} aWQ, \qquad \text{[Lemma 2]}$$

$$VP[[\alpha K_t \alpha^{\circ-} aWQ = VP[[\alpha K_t{}^{\zeta}[\alpha^{\circ-} aWQ \qquad \text{[§2(8)]}$$

$$= VP[[\alpha K_t{}^{\zeta-} cd^i caWQ, \qquad \text{[§3(l), §3(9)]}$$

$$\mathfrak{C} : VP[[\alpha K_t{}^{\zeta-} cd^i caWQ \perp\!\!\!\perp VP[[\alpha K_t{}^{\zeta-} cd^i caeWQ \qquad \text{[§3(4}_3\text{)]}$$

$$VP[[\alpha K_t{}^{\zeta-} cd^i caeWQ = VP[[\alpha K_t \alpha^{\circ-} aeWQ \qquad \text{[§3(9), §3(l), §2(8)]}$$

$$= VPcd^i c[[K_t \alpha^{\zeta-} aeWQ, \qquad \text{[§2(9), §3(l), §3(9)]}$$

$$\mathfrak{C} : VPcd^i c[[K_t \alpha^{\zeta-} aeWQ \perp\!\!\!\perp VPcd^i ca[[K_t \alpha^{\zeta-} eWQ \qquad \text{[Lemma 2]}$$

$$\perp\!\!\!\perp VPcd^i cae[[K_t \alpha^{\zeta-} [K_t \alpha^{\zeta} WQ \qquad \text{[Lemma 3]}$$

$$\perp VPcd^i ca[[K_t \alpha^{\zeta-} [K_t \alpha^{\zeta} WQ \qquad \text{[§3(4}_3\text{)]}$$

$$\perp\!\!\!\perp VPcd^i c[[K_t \alpha^{\zeta-} a[K_t \alpha^{\zeta} WQ, \qquad \text{[Lemma 2]}$$

$$VPcd^i c[[K_t \alpha^{\zeta-} a[K_t \alpha^{\zeta} WQ = VP[[\alpha K_t \alpha^{\circ-} a[K_t \alpha^{\zeta} WQ \qquad \text{[§3(l), §3(9), §2(9)]}$$

$$= VP[[\alpha K_t \alpha^{\circ-} [K_t \alpha^{\circ} WQ, \qquad \text{[§2(6)]}$$

$$\mathfrak{C} : VP[[\alpha K_t \alpha^{\circ-} [K_t \alpha^{\circ} WQ \perp\!\!\!\perp V[[\alpha K_t \alpha^{\circ-} P[K_t \alpha^{\circ} WQ \qquad \text{[Lemma 2]}$$

$$\perp\!\!\!\perp V[[\alpha K_t \alpha^{\circ-} WP[K_t \alpha^{\circ} Q, \qquad \text{[Lemma 2]}$$

$$V[[\alpha K_t \alpha^{\circ-} WP[K_t \alpha^{\circ} Q = \mathbf{S}P[K_t \alpha^{\circ} Q.$$

Combining everything, we get $\mathfrak{C} : \mathbf{S}PaQ \perp\!\!\!\perp \mathbf{S}P[K_t \alpha^{\circ} Q$, which is what was to be shown. $\qquad \square$

**Lemma 5.** *If $P$ and $Q$ are words over $A_0$, and $\mathfrak{B}^{\circ} : P \perp\!\!\!\perp Q$, then $\mathfrak{C} : \mathbf{S}P \perp\!\!\!\perp \mathbf{S}Q$.*

*Proof.* We will use induction on the number of entries in a $\mathfrak{B}^{\circ}$-sequence connecting $P$ with $Q$. If in this $\mathfrak{B}^{\circ}$-sequence only 1 entry appears, then $P = Q$, and the desired statement holds trivially. Suppose that we have proved the claim for all $\mathfrak{B}^{\circ}$-sequences containing $n$ entries (with $n \geq 1$), and suppose that a $\mathfrak{B}^{\circ}$-sequence connecting $P$ and $Q$ has $n + 1$ entries. If we denote the penultimate entry of this $\mathfrak{B}^{\circ}$-sequence by $Q_1$, then $\mathfrak{B}^{\circ} : Q_1 \perp Q$ and there exists a $\mathfrak{B}^{\circ}$-sequence with $n$ entries, which connects $P$ with $Q_1$. By the inductive hypothesis, we have $\mathfrak{C} : \mathbf{S}P \perp\!\!\!\perp \mathbf{S}Q$. But from this, as $\mathfrak{B}^{\circ} : Q_1 \perp Q$, it follows that there exists some words $M$ and $N$ over $A_0$, and an index $t$ (with $0 \leq t \leq m$), such that either $Q_1 = MaN$ and $Q = M[K_t \alpha^{\circ} N$, or else $Q_1 = M[K_t \alpha^{\circ} N$ and $Q = MaN$ by [§3(10)]. But, by Lemma 4, in both cases we have $\mathfrak{C} : \mathbf{S}Q_1 \perp\!\!\!\perp \mathbf{S}Q$, and from this it follows that also $\mathfrak{C} : \mathbf{S}P \perp\!\!\!\perp \mathbf{S}Q$. The lemma is proved. $\qquad \square$

<div align="center">

**§5**

</div>

We will define an operation which, when applied to any word $P$ over $A_1$, will produce a word $[P^\sim$, which we will call the *shape* of $P$. This operation will be defined by using right induction on $P$ over $A_1$ via the following rules:

(1) $$[\Lambda^\sim = \Lambda,$$

(2) $$[P\xi^\sim = [P^\sim[\xi^{A_o} \quad (\xi \in A_0 \cup \overline{A}_0),$$

(3) $$[Pe^\sim = [P^\sim b^i[[P^{\overline{A}_0}{}^\wedge. \qquad \text{[§3(2)]}$$

We will now establish a connection between equivalence of words by $\mathfrak{C}$ and the equivalence of words by $\mathfrak{B}^\circ$.

It is not difficult to prove the following statements:

**Lemma 1.** *If $P$ is a word over $A_1$, then $[P^\sim$ is a word over $A_0$.*

**Lemma 2.** *If $P$ is a word over $A_0 \cup \overline{A}_0$, then*

(4) $$[P^\sim = [P^{A_o}.$$

**Lemma 3.** *If $P$ is a word over $\overline{A}_0$, then*

(5) $$[Pe^\sim = b^i[P^\wedge$$

*Proof.* Indeed, we have

$$[Pe^\sim = [P^\sim b^i[[P^{\overline{A}_0}{}^\wedge \qquad\qquad \text{[(3)]}$$

$$= [P^{A_o}b^i[[P^{\overline{A}_0}{}^\wedge \qquad\qquad \text{[(4)]}$$

$$= b^i[P^\wedge \qquad\qquad \text{(as $P$ is over $\overline{A}_0$).}$$

$$\square$$

**Lemma 4.** *For any words $P$ and $Q$ over $A_1$, we have*

(6) $$[PQ^\sim = [P^\sim[[P^{\overline{A}_0}Q^\sim.$$

*Proof.* We will use right induction on $Q$ over $A_1$. For $Q = \Lambda$, equation [(6)] holds trivially, as by [(4)] we have $[[P^{\overline{A}_0}{}^\sim = [[P^{\overline{A}_0}{}^{A_o} = \Lambda$. Suppose now that for some word $Q$ over $A_1$ the equation [(6)] is true (for any $P$ over $A_1$). Let $\xi \in A_1$. We will prove, that equation [(6)] holds also when we replace $Q$ with $Q\xi$. Either $\xi \in A_0 \cup \overline{A}_0$, or $\xi = e$. Suppose first that $\xi \in A_0 \cup \overline{A}_0$. Then:

$$[PQ\xi^\sim = [PQ^\sim[\xi^{A_o} \qquad\qquad \text{[(2)]}$$

$$= [P^\sim[[P^{\overline{A}_0}Q^\sim[\xi^{A_o} \qquad\qquad \text{[(6)]}$$

$$= [P^\sim[[P^{\overline{A}_0}Q\xi^\sim, \qquad\qquad \text{[(2)]}$$

i.e. $[PQ\xi^\sim = [P^\sim[[P^{\overline{A}_0}Q\xi^{\widetilde{\phantom{x}}}$. Suppose now that $\xi = e$. Then:

$$[PQe^\sim = [PQ^\sim b^i[[PQ^{\overline{A}_0}{}^\wedge \qquad\qquad [(3)$$

$$= [P^\sim[[P^{\overline{A}_0}Q^\sim b^i[[PQ^{\overline{A}_0}{}^\wedge \qquad\qquad [(6)$$

$$= [P^\sim[[P^{\overline{A}_0}Q^\sim b^i[P^{\overline{A}_0}[[Q^{\overline{A}_0}{}^\wedge$$

$$= [P^\sim[[P^{\overline{A}_0}Q^\sim b^i[[[P^{\overline{A}_0}{}^{\overline{A}_0}[Q^{\overline{A}_0}{}^\wedge$$

$$= [P^\sim[[P^{\overline{A}_0}Q^\sim b^i[[[P^{\overline{A}_0}Q^{\overline{A}_0}{}^\wedge$$

$$= [P^\sim[[P^{\overline{A}_0}Qe^{\widetilde{\phantom{x}}} \qquad\qquad [(3)$$

and so $[PQ\xi^\sim = [P^\sim[[P^{\overline{A}_0}Q\xi^{\widetilde{\phantom{x}}}$. The lemma is proved. $\qquad\square$

Let $P, Q$, and $R$ be any words over $A_1$. We have:

$$[PQR^\sim = [PQ^\sim[[PQ^{\overline{A}_0}R^{\widetilde{\phantom{x}}} \qquad\qquad [(6)$$

$$= [P^\sim[[P^{\overline{A}_0}Q^\sim[[PQ^{\overline{A}_0}R^{\widetilde{\phantom{x}}}, \qquad\qquad [(6)$$

$$(7) \qquad\qquad [PQR^\sim = [P^\sim[[P^{\overline{A}_0}Q^\sim[[PQ^{\overline{A}_0}R^{\widetilde{\phantom{x}}}.$$

**Lemma 5.** *If $P$ and $Q$ are words over $A_1$, $\xi \in A_0$, and $\eta \in \overline{A}_0$, then $[P\xi\eta Q^\sim = [P\eta\xi Q^\sim$.*

*Proof.* Indeed, we have

$$[P\xi\eta Q^\sim = [P^\sim[[P^{\overline{A}_0}\xi\eta^{\widetilde{\phantom{x}}}[[P\xi\eta^{\overline{A}_0}Q^{\widetilde{\phantom{x}}} \qquad\qquad [(7)$$

$$= [P^\sim[[P^{\overline{A}_0}\xi\eta^{A_o}[[P\xi\eta^{\overline{A}_0}Q^{\widetilde{\phantom{x}}} \qquad\qquad [(4)$$

$$= [P^\sim\xi[[P^{\overline{A}_0}\eta Q^{\widetilde{\phantom{x}}}$$

and analogously $[P\eta\xi Q^\sim = [P^\sim\xi[[P^{\overline{A}_0}\eta Q^{\widetilde{\phantom{x}}}$. $\qquad\square$

**Lemma 6.** *If $P$ and $Q$ are words over $A_1$ and $\xi \in A_0$, then $[Pe\bar\xi\xi Q^\sim = [P\bar\xi e Q^\sim$.*

*Proof.* Indeed, we have:

$$[Pe\bar\xi\xi Q^\sim = [P^\sim[[P^{\overline{A}_0}e\bar\xi\xi^{\widetilde{\phantom{x}}}[[Pe\bar\xi\xi^{\overline{A}_0}Q^{\widetilde{\phantom{x}}} \qquad\qquad [(7)$$

$$= [P^\sim[[P^{\overline{A}_0}e\bar\xi^{\widetilde{\phantom{x}}}\xi[[Pe\bar\xi\xi^{\overline{A}_0}Q^{\widetilde{\phantom{x}}} \qquad\qquad [(2)$$

$$= [P^\sim[[P^{\overline{A}_0}e^{\widetilde{\phantom{x}}}\xi[[Pe\bar\xi\xi^{\overline{A}_0}Q^{\widetilde{\phantom{x}}} \qquad\qquad [(2)$$

$$= [P^\sim b^i[[P^{\overline{A}_0}{}^\wedge\xi[[Pe\bar\xi\xi^{\overline{A}_0}Q^{\widetilde{\phantom{x}}} \qquad\qquad [(5)$$

$$= [P^\sim b^i[[P^{\overline{A}_0}{}^\wedge\xi[[P^{\overline{A}_0}\bar\xi Q^{\widetilde{\phantom{x}}};$$

$$[P\bar\xi e Q^\sim = [P^\sim[[P^{\overline{A}_0}\bar\xi e^{\widetilde{\phantom{x}}}[[P\bar\xi e^{\overline{A}_0}Q^{\widetilde{\phantom{x}}} \qquad\qquad [(7)$$

$$= [P^\sim b^i[[P^{\overline{A}_0}\bar\xi^{\widetilde{\phantom{x}}\wedge}[[P\bar\xi e^{\overline{A}_0}Q^{\widetilde{\phantom{x}}} \qquad\qquad [(5)$$

$$= [P^\sim b^i[[P^{\overline{A}_0}{}^\wedge\xi[[P^{\overline{A}_0}\bar\xi Q^{\widetilde{\phantom{x}}}.$$

$$\qquad\square$$

**Lemma 7.** *If $P$ and $Q$ are words over $A_1$ and $\mathfrak{C} : P \perp Q$, then $[P^{\overline{A}_0} = [Q^{\overline{A}_0}$.*

*Proof.* This follows from the fact that for every defining relation of the system $\mathfrak{C}$ (see §3(4)) we have that the projection of the left-hand side onto $\overline{A}_0$ is equal to the projection on the right-hand side onto $\overline{A}_0$ (indeed, we have $[\xi\eta^{\overline{A}_0} = [\eta\xi^{\overline{A}_0} = \eta$ [§3($4_1$)], and $[e\bar{\xi}\xi^{\overline{A}_0} = [\bar{\xi}e^{\overline{A}_0} = \bar{\xi}$ [§3($4_2$)], and finally $[cd^ica^{\overline{A}_0} = [cd^icae^{\overline{A}_0} = cd^ic$ [§3($4_3$)]).    □

**Lemma 8.** *If $P$ is a word over $A_0$, $R$ is a word over $A_1$, and $\mathfrak{C} : \mathbf{S}P \perp\!\!\!\perp R$, then $[R^{\overline{A}_0} = \mathbf{S}$.*

The proof of this is using Lemma 7 and induction on the number of entries in a $\mathfrak{C}$-sequence connecting $\mathbf{S}P$ with $R$.

**Lemma 9.** *If the word $Pcd^ic$ is a prefix of $\mathbf{S}$, where $P$ is a word over $\overline{A}_0$, then*

$$\mathfrak{B}^{\bigcirc} : ab^i[Pcd^ic^{\wedge} \perp\!\!\!\perp a.$$

*Proof.* There exists some word $Q$ over $\overline{A}_0$ such that $\mathbf{S} = Pcd^icQ$. Then we have

$$[\mathbf{S}^{\wedge} = [Pcd^icQ^{\wedge} = [P^{\wedge}ab^ia[Q^{\wedge} = [P^{\wedge}[\alpha^{\bigcirc}[Q^{\wedge}. \qquad\qquad [§3(9)]$$

On the other hand, by [§3(3), §3(11)], we have $[\mathbf{S}^{\wedge} = [\alpha K_0\alpha K_1\alpha\cdots K_m\alpha^{\bigcirc}$ and furthermore also $[\alpha K_0\alpha K_1\alpha\cdots K_m\alpha^{\bigcirc} = [P^{\wedge}[\alpha^{\bigcirc}[Q^{\wedge}$. By [§2, Lemma 4], there exist words $T$ and $U$ over $B$, such that

$$[T^{\varsigma} = [P^{\wedge}, [U^{\beth} = [Q^{\wedge}, \text{ and } T\alpha U = \alpha K_0\alpha K_1\alpha\cdots K_m\alpha.$$

But $K_0, K_1, \ldots, K_m$ are words over $A$, and $\alpha \notin A$, so $T$ must be one of the words $\Lambda, \alpha K_0, \alpha K_0\alpha K_1, \ldots$, and $\alpha K_0\alpha K_1\cdots\alpha K_m$, and consequently by [§3(7), §3(8)], we have $\mathfrak{B} : T \perp\!\!\!\perp \Lambda$. Then by [§3(7)] it follows that $\mathfrak{B} : \alpha T\alpha \perp\!\!\!\perp \Lambda$, and by [§3, Lemma 5] also $\mathfrak{B}^{\bigcirc} : [\alpha T\alpha^{\bigcirc} \perp\!\!\!\perp [\Lambda^{\bigcirc}$. But $[\Lambda^{\bigcirc} = a$, and

$$[\alpha T\alpha^{\bigcirc} = [\alpha^{\varsigma}[T\alpha^{\bigcirc} = [\alpha^{\varsigma}[T^{\varsigma}[\alpha^{\bigcirc} = ab^i[P^{\wedge}ab^ia = ab^i[Pcd^ic^{\wedge}; \qquad [§2(8), §3(9)]$$

and hence $\mathfrak{B}^{\bigcirc} : ab^i[Pcd^ic^{\wedge} \perp\!\!\!\perp a$, which is what was to be shown.    □

**Lemma 10.** *If $P$ and $Q$ are words over $A_1$, and $[Pcd^icaQ^{\overline{A}_0} = \mathbf{S}$, then*

$$\mathfrak{B}^{\bigcirc} : [Pcd^icaQ^{\sim} \perp\!\!\!\perp [Pcd^icaeQ^{\sim}.$$

*Proof.* We have:

$$[Pcd^icaQ^{\sim} = [P^{\sim}[[P^{\overline{A}_0}cd^ica^{\sim}[[Pcd^ica^{\overline{A}_0}Q^{\sim} \qquad\qquad [(7)]$$

$$= [P^{\sim}[[P^{\overline{A}_0}cd^ica^{A_o}[[Pcd^ica^{\overline{A}_0}Q^{\sim} \qquad\qquad [(4)]$$

$$= [P^{\sim}a[[P^{\overline{A}}cd^icQ^{\sim};$$

$$[Pcd^icaeQ^{\sim} = [P^{\sim}[[P^{\overline{A}_0}cd^icae^{\sim}[[Pcd^icae^{\overline{A}_0}Q^{\sim} \qquad\qquad [(7)]$$

$$= [P^{\sim}[[P^{\overline{A}_0}cd^ica^{\sim}b^i[[[P^{\overline{A}_0}cd^ica^{\overline{A}_0}{}^{\wedge}[[Pcd^icae^{\overline{A}_0}Q^{\sim} \qquad [(3)]$$

$$= [P^{\sim}[[P^{\overline{A}_0}cd^ica^{\sim}b^i[[P^{\overline{A}_0}cd^ic^{\wedge}[[P^{\overline{A}_0}cd^icQ^{\sim}$$

$$= [P^{\sim}[[P^{\overline{A}_0}cd^ica^{A_o}b^i[[P^{\overline{A}_0}cd^ic^{\wedge}[[P^{\overline{A}_0}cd^icQ^{\sim} \qquad\qquad [(4)]$$

$$= [P^{\sim}ab^i[[P^{\overline{A}_0}cd^ic^{\wedge}[[P^{\overline{A}_0}cd^icQ^{\sim}.$$

Hence, the word $[Pcd^{i}caeQ^{\sim}$ is obtained from $[Pcd^{i}caQ^{\sim}$ by replacing some occurrence of the letter $a$ by the word $ab^{i}[[P^{\overline{A_0}}cd^{i}c^{\wedge}$. But $\mathbf{S} = [Pcd^{i}caQ^{\overline{A_0}} = [P^{\overline{A_0}}cd^{i}c[Q^{\overline{A_0}},$ the word $[P^{\overline{A_0}}cd^{i}c$ is a prefix of $\mathbf{S}$, and, by Lemma 9, we have $\mathfrak{B}^{\circ} : ab^{i}[[P^{\overline{A_0}}cd^{i}c^{\wedge} \perp\!\!\!\perp a.$ From this, it follows that $\mathfrak{B}^{\circ} : [Pcd^{i}caQ^{\sim} \perp\!\!\!\perp [Pcd^{i}caeQ^{\sim}.$ □

**Lemma 11.** *If $P$ is a word over $A_0$, and $\mathfrak{C} : \mathbf{S}P \perp\!\!\!\perp Q$, then $\mathfrak{B}^{\circ} : [Q^{\sim} \perp\!\!\!\perp P$.*

*Proof.* We will use induction on the number of entries in a $\mathfrak{C}$-sequence connecting $\mathbf{S}P$ with $Q$. If there is only one entry in this sequence, then $\mathbf{S}P = Q$, and so by [(4)] we have $[Q^{\sim} = [Q^{A_0} = [\mathbf{S}P^{A_0} = P$. Suppose that our claim is proved in the cases that the $\mathfrak{C}$-sequence connecting $\mathbf{S}P$ with $Q$ has $n$ entries (with $n \geq 1$); suppose now that the $\mathfrak{C}$-sequence connecting $\mathbf{S}P$ with $Q$ has $n+1$ entries. If we denote the penultimate term of this $\mathfrak{C}$-sequence by $Q_1$, then $\mathfrak{C} : Q_1 \perp Q$ and $\mathfrak{S}P$ is connected with $Q_1$ by a $\mathfrak{C}$-sequence with $n$ entries. By the inductive hypothesis, we have $\mathfrak{B}^{\circ} : [Q_1^{\sim} \perp\!\!\!\perp P$. From the fact that $\mathfrak{C} : Q_1 \perp Q$, it follows that there exists some words $T$ and $U$ over $A_1$ such that one of the following holds:

(a) one of $Q_1$ and $Q$ equals $T\xi\eta U$, the other $T\eta\xi U$, where $\xi \in A_0, \eta \in \overline{A}_0$;  [§3($4_1$)]
(b) one of $Q_1$ and $Q$ equals $Te\xi\xi U$, the other $T\bar{\xi}eU$, where $\xi \in A_0$;     [§3($4_2$)]
(c) one of $Q_1$ and $Q$ equals $Tcd^{i}caU$, the other $Tcd^{i}caeU$.         [§3($4_3$)]

In case (a) and case (b), using Lemma 5 and Lemma 6, respectively, we have that $[Q^{\sim} = [Q_1^{\sim}$, and consequently that $\mathfrak{B}^{\circ} : [Q^{\sim} \perp\!\!\!\perp P$. In case (c), then by using the fact that $\mathfrak{C} : \mathbf{S} \perp\!\!\!\perp Q$ and $\mathfrak{C} : \mathbf{S}P \perp\!\!\!\perp Q_1$, it follows that $\mathfrak{C} : \mathbf{S}P \perp\!\!\!\perp Tcd^{i}caU$ and, by Lemma 8, also $[Tcd^{i}caU^{\overline{A_0}} = \mathbf{S}$. Then by Lemma 10, we have $\mathfrak{B}^{\circ} : [Tcd^{i}caU^{\sim} \perp\!\!\!\perp [Tcd^{i}caeU^{\sim},$ i.e. $\mathfrak{B}^{\circ} : [Q^{\sim} \perp\!\!\!\perp [Q_1^{\sim}$, and then again we have $\mathfrak{B}^{\circ} : [Q^{\sim} \perp\!\!\!\perp P$. Our lemma is proved. □

**Lemma 12.** *If $P$ and $Q$ are words over $A_0$, and $\mathfrak{C} : \mathbf{S}P \perp\!\!\!\perp \mathbf{S}Q$, then $\mathfrak{B}^{\circ} : P \perp\!\!\!\perp Q$.*

Indeed, by Lemma 11, we have $\mathfrak{B}^{\circ} : [\mathbf{S}Q^{\sim} \perp\!\!\!\perp P$, so by [(4)] also $[\mathbf{S}Q^{\sim} = [\mathbf{S}Q^{\overline{A_0}} = Q$ and $\mathfrak{B}^{\circ} : P \perp\!\!\!\perp Q$.

## §6

We can now prove the unsolvability of the equivalence problem in the calculus $\mathfrak{C}$.

From [§4, Lemma 5] and [§5, Lemma 12] the following statements clearly follow (already stated at the end of §3):

**Lemma 1.** *If $P$ and $Q$ are words over $A_0$, then $\mathfrak{B}^{\circ} : P \perp\!\!\!\perp Q$ if and only if $\mathfrak{C} : \mathbf{S}P \perp\!\!\!\perp \mathbf{S}Q$.*

Furthermore, from [§3, Lemma 2] and [§3, Lemma 4] follows:

**Lemma 2.** *If $P$ and $Q$ are words over $A$, then $\mathfrak{A} : P \perp\!\!\!\perp Q$ if and only if $\mathfrak{B} : P \perp\!\!\!\perp Q$.*

Finally, from [§3, Lemma 5] and [§3, Lemma 7] we have:

**Lemma 3.** *If $P$ and $Q$ are words over $A$, then $\mathfrak{B} : P \perp\!\!\!\perp Q$ if and only if $\mathfrak{B}^{\circ} : [P^{\circ} \perp\!\!\!\perp [Q^{\circ}$.*

Thus, combining Lemmas 1,2,3 we find:

**Lemma 4.** *If $P$ and $Q$ are words over $A$, then $\mathfrak{A} : P \perp\!\!\!\perp Q$ if and only if $\mathfrak{C} : \mathbf{S}[P^{\circ} \perp\!\!\!\perp \mathbf{S}[Q^{\circ}$.*

Remembering that $\mathfrak{A}$ can be taken as an arbitrary special associative calculus (see §3), we thus see by Lemma 4 that *the equivalence problem for any pair of words over any special associative calculus can be reduced to the equivalence problem in $\mathfrak{C}$ of some two words,*

*which can be effectively constructed from the given words and the defining relations of the given calculus.*

Let us call the *uniform equivalence problem for special associative calculi* the problem of constructing an algorithm (§1) which could take as input any special associative calculus and any pair of words over the alphabet of this calculus and decide whether or not these words are equivalence in the calculus or not. In view of the above, the uniform equivalence problem for special associative calculi can be reduced to the equivalence problem in $\mathfrak{C}$. Consequently, we have the following:

**Lemma 5.** *If the uniform equivalence problem for special associative calculi is unsolvable, then the equivalence problem in $\mathfrak{C}$ is unsolvable.*

P. S. Novikov [4] has constructed a group with a finite number of generators and defining relations, in which the identity problem is unsolvable. Denote the alphabet which consists of all the generators of this group and their inverses by $A$; and denote by $\xi^{-1}$ the inverse element of the generator $\xi$ (for $\xi \in A$); if $P$ is a word over $A$, and $P = \xi_1 \xi_2 \cdots \xi_t$, then we denote by $\lceil P^{-1}$ the word $\xi_t^{-1} \cdots \xi_2^{-1} \xi_1^{-1}$, i.e.

$$\lceil P^{-1} = \mathbf{S}_{\cdots \xi^{-1} \cdots}^{\cdots \xi \cdots} \lceil P^{\smile} \rfloor.$$

Clearly, in the group we have $P \lceil P^{-1} = 1$ for any word $P$ over $A$. Let the defining system of the group have the form

$$\{ M_t = N_t, \quad (1 \leq t \leq s).$$

We construct an associative calculus $\mathfrak{A}$ over the alphabet $A$ with the following defining system:

$$\begin{cases} M_t \lceil N_t^{-1} \longleftrightarrow & (1 \leq t \leq s) \\ \xi \xi^{-1} \longleftrightarrow & (\xi \in A). \end{cases}$$

It is not difficult to show that two words are equal in the group if and only if they are equivalent in $\mathfrak{A}$. Hence the problem of equivalence in $\mathfrak{A}$ is unsolvable. But $\mathfrak{A}$ is a special associative calculus, and if the uniform equivalence problem for special associative calculi were solvable then the equivalence problem in $\mathfrak{A}$ must also be solvable. Hence:

**Lemma 6.** *The uniform equivalence problem for special associative calculi is unsolvable.*

From Lemmas 5 and 6 now follow:

**Theorem 1.** *The equivalence problem in $\mathfrak{C}$ is unsolvable.*

**Corollary.** *The equivalence problem is unsolvable in the calculus $\mathfrak{C}_1$ [§1(1)].*

The proof of Theorem 1 relies on the theorem by P. S. Novikov on the existence of groups with unsolvable identity problem, which itself is proved in a very complicated way. The proof of Theorem 1 will thus be simplified if we can find a direct proof of Lemma 6. However, the author is not aware of such a proof.

## §7

We will now construct, finally, a calculus for which the problem of equivalence with a particular word is unsolvable.

Fix some non-negative integer $j$, different from $i$, and define the associative calculus $\mathfrak{C}'$ over the alphabet $A_1$ and with the following system of defining relations (written in an abbreviated way):

$$(1) \quad \begin{cases} (1_1) & \xi\eta \longleftrightarrow \eta\xi & (\xi \in A_0, \eta \in \overline{A}_0) \\ (1_2) & e\bar{\xi}\xi \longleftrightarrow \bar{\xi}e & (\xi \in A_0) \\ (1_3) & cd^i ca \longleftrightarrow cd^i cae \\ (1_4) & \eta ab^j ab^j a \longleftrightarrow ab^j ab^j a & (\eta \in \overline{A}_0). \end{cases}$$

When $i = 1, j = 0$ this calculus is simply the calculus $\mathfrak{C}_2$ (see §1(2)).

**Lemma 1.** *If $P$ and $Q$ are words over $A_1$, and $\mathfrak{C} : P \perp\!\!\!\perp Q$, then $\mathfrak{C}' : P \perp\!\!\!\perp Q$.*

This is obvious by comparing [§3(4)] with [$(1_1)$, $(1_2)$, $(1_3)$].

**Lemma 2.** *If $P$ is a word over $\overline{A}_0$, then $\mathfrak{C}' : Pab^j ab^j a \perp\!\!\!\perp ab^j ab^j a$.*

This is proved by right induction over $\overline{A}_0$ and using [$(1_4)$].

**Lemma 3.** *If a $\mathfrak{C}'$-sequence is such that no word containing an occurrence of the word $ab^j ab^j a$ appears in it, then this $\mathfrak{C}'$-sequence is a $\mathfrak{C}$-sequence.*

*Proof.* Indeed, if in this $\mathfrak{C}'$-sequence there is a pair of words which is connected by an application of [$(1_4)$], then both the words in this pair will contain an occurrence of $ab^j ab^j a$, which is impossible. Hence, in this $\mathfrak{C}'$-sequence the only relations used are those of the form [$(1_1)$, $(1_2)$, and $(1_3)$], which means that it is a $\mathfrak{C}$-sequence [§3(4)].    □

Suppose that $\mathfrak{A}'$ is any special associative calculus and let $A'$ be the alphabet of $\mathfrak{A}'$. Let $\beta$ be a letter not occurring in $A'$, and let $A = A' \cup \{\beta\}$. Let $\mathfrak{A}$ be the associative calculus over $A$ with the same defining relations as $\mathfrak{A}'$. Then $\mathfrak{A}$ is also a special associative calculus.

**Lemma 4.** *If $P$ and $Q$ are words over $A'$ and $\mathfrak{A}' : P \perp\!\!\!\perp Q$, then $\mathfrak{A} : P \perp\!\!\!\perp Q$.*

**Lemma 5.** *If $P, Q$, and $R$ are words over $A'$, $V$ is a word over $A$, and $\mathfrak{A} : P\beta Q\beta R \perp V$, then there exists some word $X$ over $A'$ such that one of the following holds:*

*(a)* $\mathfrak{A}' : P \perp X$ *and* $V = X\beta Q\beta R$,
*(b)* $\mathfrak{A}' : Q \perp X$ *and* $V = P\beta X\beta R$,
*(c)* $\mathfrak{A}' : R \perp X$ *and* $V = P\beta Q\beta X$.

The proof follows directly from [2, I.§4.5.4].

**Lemma 6.** *If $P, Q$, and $R$ are words over $A'$, $V$ is a word over $A$, and $\mathfrak{A} : P\beta Q\beta R \perp\!\!\!\perp V$, then there exists some words $P_1, Q_1$, and $R_1$ over $A'$ such that $V = P_1\beta Q_1\beta R_1$, $\mathfrak{A}' : P \perp\!\!\!\perp P_1$, $\mathfrak{A}' : Q \perp\!\!\!\perp Q_1$, and $\mathfrak{A}' : R \perp\!\!\!\perp R_1$.*

The proof is proved using Lemma 5 and induction on the number of entries in an $\mathfrak{A}'$-sequence connecting $P\beta Q\beta R$ with $V$.

We now apply the same procedure to $\mathfrak{A}$ as we did in §3: construct an alphabet $B$, define an encoding from $B$ into $A_0$, define the calculi $\mathfrak{B}$ and $\mathfrak{B}^{\circ}$ and the word $\mathbf{S}$. When doing this, when defining the encoding from $B$ into $A_0$, we will set $n(\beta) = j$ (this is possible, as $j \neq i$). Then

$$(2) \qquad\qquad [\beta\beta^{\circ} = ab^j ab^j a. \qquad\qquad [\S2(1), \S2(6)]$$

**Lemma 7.** *If $P$ is a word over $A'$ and $\mathfrak{A}' : P \perp\!\!\!\perp \Lambda$, then $\mathfrak{C}' : \mathbf{S}[\beta P \beta^{\bigcirc} \perp\!\!\!\perp ab^j ab^j a$.*

*Proof.* By Lemma 4, we have $\mathfrak{A} : P \perp\!\!\!\perp \Lambda$, hence $\mathfrak{A} : \beta P \beta \perp\!\!\!\perp \beta\beta$. By [§6, Lemma 4] we have $\mathfrak{C} : \mathbf{S}[\beta P \beta^{\bigcirc} \perp\!\!\!\perp \mathbf{S}[\beta\beta^{\bigcirc}$. Then by Lemma 1, we have $\mathfrak{C}' : \mathbf{S}[\beta P \beta^{\bigcirc} \perp\!\!\!\perp \mathbf{S}[\beta\beta^{\bigcirc}$, i.e. by [(2)] we have $\mathfrak{C}' : \mathbf{S}[\beta P \beta^{\bigcirc} \perp\!\!\!\perp \mathbf{S}ab^j ab^j a$. But $\mathbf{S}$ is a word over $\overline{A}_0$, so by Lemma 2 it follows that $\mathfrak{C}' : \mathbf{S}ab^j ab^j a \perp\!\!\!\perp ab^j ab^j a$ and hence also $\mathfrak{C}' : \mathbf{S}[\beta P \beta^{\bigcirc} \perp\!\!\!\perp ab^j ab^j a$, which is what was to be shown. □

**Lemma 8.** *If $P$ is a word over $A'$ and $\mathfrak{C}' : \mathbf{S}[\beta P \beta^{\bigcirc} \perp\!\!\!\perp ab^j ab^j a$, then $\mathfrak{A}' : P \perp\!\!\!\perp \Lambda$.*

*Proof.* In a $\mathfrak{C}'$-sequence connecting $\mathbf{S}[\beta P \beta^{\bigcirc}$ with $ab^j ab^j a$ there are entries containing occurrences of $ab^j ab^j a$, including e.g. the last entry. Suppose that $Mab^j ab^j aN$ is the first of the entries in this $\mathfrak{C}'$-sequence which contains an occurrence of $ab^j ab^j a$. Then the part of this $\mathfrak{C}'$-sequence running from the beginning to the word $Mab^j ab^j aN$ contains only one entry containing an occurrence of $ab^j ab^j a$. By Lemma 3, this part of the $\mathfrak{C}'$-sequence will be a $\mathfrak{C}$-sequence, and $\mathfrak{C} : \mathbf{S}[\beta P \beta^{\bigcirc} \perp\!\!\!\perp Mab^j ab^j aN$. By [§5, Lemma 1] we have $\mathfrak{B}^{\bigcirc} : [Mab^j ab^j aN^{\sim} \perp\!\!\!\perp [\beta P \beta^{\bigcirc}$. But:

$$[Mab^j ab^j aN^{\sim} = [M^{\sim}[[M^{\overline{A}_0}ab^j ab^j a^{\sim}[[Mab^j ab^j a^{\overline{A}_0}N^{\sim} \qquad \text{[§5(7)]}$$

$$= [M^{\sim}[[M^{\overline{A}_0}ab^j ab^j a^{A_o}[[Mab^j ab^j a^{\overline{A}_0}N^{\sim} \qquad \text{[§5(4)]}$$

$$= [M^{\sim}ab^j ab^j a[[M^{\overline{A}_0}N^{\sim}$$

$$= [M^{\sim}[\beta\beta^{\bigcirc}[[M^{\overline{A}_0}N^{\sim} \qquad \text{[(2)]}$$

and $\mathfrak{B}^{\bigcirc} : [\beta P \beta^{\bigcirc} \perp\!\!\!\perp [M^{\sim}[\beta\beta^{\bigcirc}[[M^{\overline{A}_0}N^{\sim}$. In view of [§3, Lemma 6] there exist a word $Q$ over $B$ such that $\mathfrak{B} : \beta P \beta \perp\!\!\!\perp Q$ and $[M^{\sim}[\beta\beta^{\bigcirc}[[M^{\overline{A}_0}N^{\sim} = [Q^{\bigcirc}$. Then from [§2, Lemma 4] it follows that there exists some words $T$ and $U$ over $B$ such that $[T^{\complement} = [M^{\sim}, [U^{\partial} = [[M^{\overline{A}_0}N^{\sim}$, and $Q = T\beta\beta U$. Hence $\mathfrak{B} : \beta P \beta \perp\!\!\!\perp T\beta\beta U$, and by [§3, Lemma 4] also $\mathfrak{A} : [\beta P \beta^A \perp\!\!\!\perp [T\beta\beta U^A$. Since $P$ is a word over $A'$, we have $\mathfrak{A} : \beta P \beta \perp\!\!\!\perp [T^A\beta\beta[U^A$, and by Lemma 6 there are some words $X, Y, Z$ over $A'$ such that $[T^A\beta\beta[U^A = X\beta Y\beta Z$, $\mathfrak{A}' : \Lambda \perp\!\!\!\perp X$, $\mathfrak{A}' : P \perp\!\!\!\perp Y$, and $\mathfrak{A}' : \Lambda \perp\!\!\!\perp Z$. The word $X\beta Y\beta Z$ contains only two occurrences of $\beta$, as $[T^A$ and $[U^A$ do not contain occurrences of $\beta$, and $[T^A = X, \Lambda = Y$, and $[U^A = Z$. From this, we find that $Y = \Lambda$ and $\mathfrak{A}' : P \perp\!\!\!\perp Y$, from which we also conclude that $\mathfrak{A}' : P \perp\!\!\!\perp \Lambda$. This proves the lemma. □

From Lemma 7 and Lemma 8 we now have:

**Lemma 9.** *If $P$ is a word over $A'$, then $\mathfrak{A}' : P \perp\!\!\!\perp \Lambda$ if and only if $\mathfrak{C}' : \mathbf{S}[\beta P \beta^{\bigcirc} \perp\!\!\!\perp ab^j ab^j a$.*

Analogously to the uniform equivalence problem for special associative calculi, we may formulate the ***uniform problem of equivalence with the empty word for special associative calculi*** as the problem of finding an algorithm which can, when given any special associative calculus and any word over the alphabet of this calculus, decide whether or not this word is equivalent to the empty word in the given calculus. By Lemma 9, it follows that this problem can be reduced to the problem of equivalence with the word $ab^j ab^j a$ in the calculus $\mathfrak{C}'$. Thus we have the following:

**Lemma 10.** *If the uniform problem of equivalence with the empty word for special associative calculi is unsolvable, then the problem of equivalence with the word $ab^j ab^j a$ is unsolvable in $\mathfrak{C}'$.*

Following an argument similar to that in the preceding section, it follows that a special instance of the uniform problem of equivalence with the empty word for special associative calculi is the problem of being equal to the identity element in the group of P. S. Novikov [4]. But obviously, the question of whether two words $M$ and $N$ are equal in a group is equivalent to the question of whether the word $M[N^{-1}$ is equal to the identity element. Hence the problem of being equal to the identity element in this group is also unsolvable. Consequently, we have the following:

**Lemma 11**. *The uniform problem of equivalence with the empty word for special associative calculi is unsolvable.*

Finally, from Lemma 10 and Lemma 11 we find:

**Theorem 2**. *In the calculus $\mathfrak{C}'$ the problem of equivalence with the word $ab^{j}ab^{j}a$ is unsolvable.*

**Corollary 1**. *In the calculus $\mathfrak{C}_2$ (defined in §1(2)) the problem of equivalence with the word $aaa$ is unsolvable.*


 
 \end{document}